%% file: Paper.tex
\documentclass[11pt,reqno]{amsart}

\usepackage[T1]{fontenc}
\usepackage[english]{babel}
\usepackage{mathtools,amssymb,amsthm}
\usepackage{graphicx}
\usepackage{tensor}
\usepackage{accents}
\usepackage{enumitem}
\usepackage{fullpage} % Sets all margins to 1 in.  
\usepackage{microtype}
\usepackage[colorlinks=true,linkcolor=blue,citecolor=blue,urlcolor=blue]{hyperref}
\usepackage[nameinlink,noabbrev]{cleveref}

\newcommand{\lmd}{\lambda}
\newcommand{\ls}{\lesssim}

\newcommand{\ep}{\varepsilon}

\newcommand{\dd}{\,\mathrm{d}}

\newcommand{\vn}[1]{\|#1\|}
\newcommand{\vm}[1]{\left|#1\right|}
\newcommand{\mb}{\mathbb}
\newcommand*\jb[1]{\left\langle #1 \right\rangle} % Japanese Bracket
\newcommand{\pt}{\partial}
\newcommand{\Hn}{\mathbb{H}^n}

\newtheorem{theorem}{Theorem}[section]
\newtheorem{corollary}[theorem]{Corollary}
\newtheorem{lemma}[theorem]{Lemma}
\newtheorem{proposition}[theorem]{Proposition}
\theoremstyle{definition}
\newtheorem{definition}[theorem]{Definition}

\newtheorem{remark}[theorem]{Remark}
\numberwithin{equation}{section}

\title{Bilinear Strichartz estimates on the hyperbolic space and the cubic NLS}
\author{Cristian Gavrus}
\email{gcd2006@gmail.com}
\begin{document}
\maketitle

\begin{abstract}
We prove bilinear Strichartz estimates for the Schr\"odinger equation on the hyperbolic space $\Hn$.
These bilinear estimates are then applied, as part of an $I$-method argument, to establish global well-posedness for
the cubic defocusing nonlinear Schr\"odinger equation on the hyperbolic plane for regularities above $ 4/7 $.
The $I$-method is implemented in physical space using a heat-based modified energy. 
\end{abstract}

\setcounter{tocdepth}{1}
\tableofcontents

\input{sections/introduction}

\input{sections/proof_bilinear}

\input{sections/ac_main_thm_proof.tex}

\input{sections/proof_ac_law.tex}

\input{sections/proof_lwp.tex}

\input{sections/appendix}

\bibliographystyle{abbrv}
\bibliography{references}

\end{document}

%% file: sections/introduction.tex
\addtocontents{toc}{\protect\setcounter{tocdepth}{2}}
\section{Introduction}

Consider the linear Schr\"odinger equation 
\begin{equation}\label{eq:schrodinger}
 i \partial_t u + \Delta u = 0, \qquad \qquad x \in \mathcal{M}
\end{equation}
and the cubic defocusing nonlinear Schr\"odinger equation
\begin{equation} \label{eq:nls}
i\partial_t u + \Delta u = |u|^2 u , \qquad \quad x \in \mathcal{M} %, \qquad \qquad x \in \mathcal{M} % \quad \text{on } \mathbb{H}^2
\end{equation}
on a manifold $ \mathcal{M} $. 

The first purpose of this paper is to prove a bilinear Strichartz estimate for \eqref{eq:schrodinger} 
on the hyperbolic space $ \mathcal{M}=\mathbb{H}^n $, 
\begin{equation}\label{eq:BS}
\left\|u v \right\|_{L^2_{t,x}(\mathbb{R}\times\mathcal{M})}
\lesssim
M^{\frac{n-1}{2}}\,N^{-\frac12}\,
\|f\|_{L^2_x}\,\| g\|_{L^2_x}
\end{equation}
for solutions $ u, v $ localized at dyadic frequencies $ N \geq M \geq 1 $, with initial data $ f, g $. 

The second goal is to show that this bilinear gain is strong enough to run an $I$-method argument for the cubic NLS \eqref{eq:nls} 
on $ \mathbb{H}^2 $ provided the euclidean Fourier modified energy is replaced by a more geometric heat-flow almost-conserved energy.

A common theme in both arguments is bypassing the absence of the euclidean identity relating multiplication in physical space and convolution in Fourier space.

\

\subsection{The Bilinear Strichartz estimate on $ \mathbb{H}^n $}

On $ \mathcal{M}=\mathbb{R}^n $ the bound \eqref{eq:BS} is known as Bourgain's bilinear refinement of Strichartz estimates.
Introduced in \cite{bourgain1998refinements}, this estimate has inspired and influenced a wide range of nonlinear developments including:
rough NLS solutions \cite{colliander2002almost, colliander2001bilinear, colliander2008global}, % (Fourier high-low truncation \cite{}, I-method \cite{}), energy-critical NLS \cite{} (see also \cite{}, \cite{}), 
transversality \cite{tao2001multilinear}, %(Tao)
lateral smoothing  \cite{ionescu2006low},  % on hyperplanes
geometric and coupled equations \cite{bejenaru20092d,liu2014local,bejenaru2011global}, % \cite{}, \cite{} % (Schrodinger maps, Chern–Simons–Schrodinger, Zakharov)
compact manifold extensions \cite{burq2005bilinear,jiang2011bilinear,hani2012bilinear,deng2024bilinear,deng2025sharp}, %\cite{}, \cite{}, \cite{}, % (sphere, Hani, periodic, Zoll, wave-guide etc),  % and semi-classical
bilinear virial and nonlinear formulations \cite{planchon2009bilinear, tao2010physical, dodson2019defocusing, dodson2012bilinear}, %\cite{planchon2009bilinear}, \cite{tao2010physical}, \cite{dodson2019defocusing}, \cite{dodson2012bilinear},  %(Planchon-Vega, Ben book, Ben wave maps paper),  interaction Morawetz 
%angular refinements \cite{}, 
probabilistic applications \cite{burqBilAS, benyi2015probabilistic, dodson2019almost}, % \cite{burqBilAS}, \cite{benyi2015probabilistic}, \cite{dodson2019almost}, %(Benyi, Oh, Pocovnicu; Dodson–Lührmann–Mendelson)
to name just a few. 

\

The wide applicability of the euclidean estimate and the interest in dispersive equations on $\mathbb{H}^n$ 
\cite{banica2007nonlinear,anker2009nonlinear,ionescu2009semilinear,metcalfe2012dispersive,ionescu2012global,lawrie2023localsmoothing,lawrie2023asymptotic,Wilson2023GlobalWA}
motivate the investigation of bilinear analogues on the hyperbolic space.

\begin{theorem}[Bilinear Strichartz on $\mathbb{H}^n$]\label{thm:bilinear-hyperbolic}
Suppose $ u, v : \mb{R} \times \mathbb{H}^n \to \mb{C} $ are solutions to the linear Schr\"odinger equation \eqref{eq:schrodinger} with initial data $ f, g \in L^2(\mathbb{H}^n) $ respectively, for $ n \geq 2 $.
Assume that $ u, v $ are spectrally localized at dyadic frequencies $ N, M $ respectively, with $ N \geq M\geq 1$. Then 
\begin{equation}\label{BS}
\left\|u v \right\|_{L^2_{t,x}(\mathbb{R}\times \mathbb{H}^n)}
\lesssim
M^{\frac{n-1}{2}}\,N^{-\frac12}\,
\|f\|_{L^2_x}\,\| g\|_{L^2_x}.
\end{equation}
and the same holds for $ \bar{u} v $ and $ u \bar{v} $. The estimate \eqref{BS} is sharp. 
\end{theorem}

More precisely, the spectral localization means that the Helgason-Fourier transform of $f$, $ \tilde{f}(\lambda,\omega) $ is supported in intervals such as $ \lambda \in [N/4, 4N] $ (or $[0,4]$ when $ N=1$) and similarly for $g$ and $ M $. 

When $N \simeq M$ \eqref{BS} follows from linear Strichartz estimates \cite{anker2009nonlinear,banica2007nonlinear,ionescu2009semilinear} (such as $L^4_{t,x} $ when $n=2$). 
The genuinely bilinear regime occurs % at high-low frequencies interactions 
when $ N \gg M $: % which are suppressed by a scale ratio, reflecting a transversality gain:
waves at separated frequencies interact less strongly than suggested by linear Strichartz estimates. %by a direct application of 
The $ (M/N)^{\frac12} $ gain expresses short time interactions between a high-frequency packet rapidly crossing the low frequency. 

Multiplying two waves on $ \mb{R}^n $ corresponds to convolving their spacetime Fourier supports, which lie on two portions of the paraboloid whose normals differ by $\simeq N$ when $ N \gg M $. 
The overlap geometry of the thickened paraboloids with large separation of normals is the restriction-theoretic origin of the euclidean factor $ M^{\frac{n-1}{2}} N^{-\frac12} $.
See \cite{tao2001multilinear} (Sec. 7, 11), \cite{koch2014dispersive} (Thm. 5.1, 5.6 in Ch. 5).

\

On $ \mathbb{H}^n $ it is not meaningful to view the product of two functions as some sort of convolution for their Fourier transforms.
However, it is shown by the proof of Theorem \ref{thm:bilinear-hyperbolic} that this obstruction is not essential. 
One may linearize the problem by lifting it to the tensor product space, and $TT^*$ becomes convolution with a radial kernel.
Its boundedness on $ L^2(\mb{R}\times \mathbb{H}^n) $ is reduced to estimating spectral integrals analogous to those 
arising in the euclidean “convolution of hypersurface measures" argument, concluding with the help of 
bounds for the linearization coefficients of spherical functions. %Jacobi triple product coefficients.
The initial data are not assumed radial.

\

To make the estimate \eqref{BS} suited for nonlinear applications, the following standard extensions are in order. 
First, replace the compact spectral support condition by the more geometric heat flow localization, as
$ e^{s\Delta}$ have good $L^p(\mathbb{H}^n)$ properties, unlike compactly supported spectral multipliers \cite{anker1990p,clerc1974lp,ionescu2002singular}.  
Secondly, one may extend all estimates to $ X^{s,b} $ inputs using a well-known transfer principle. 
The definition and properties of $ X^{s,b}(\mathbb{H}^n) $ spaces are presented in section \ref{sec:xsb_spaces}.

\begin{corollary}[Bilinear Strichartz - heat flow localization] \label{cor:bilinear_strichartz_variations} 
    (1) Suppose $ u, v : \mb{R} \times \mathbb{H}^n \to \mb{C} $ are solutions to
     the linear Schr\"odinger equation \eqref{eq:schrodinger} with initial data $ f, g \in L^2(\mathbb{H}^n) $ respectively. \
    Then, for $ 0<r \leq s $ and $ \theta > \frac{1}{4} $  one has
\begin{equation}
    \label{eq:heat_bil}
    \vn{(1-e^{r\Delta})u \, e^{s\Delta} v}_{L^2_{t,x}} + \vn{(-r \Delta)^{\theta} e^{r\Delta} u \, e^{s\Delta} v}_{L^2_{t,x}} \ls r^{\frac14} s^{-\frac{n-1}{4}}  \|f\|_{L^2} \|g\|_{L^2}. 
\end{equation}
  (2)   (Inhomogeneous estimate) All bounds extend to $ X^{0,b} $ inputs for $ b > \frac12 $. In particular, for $r \leq s$ 
    \begin{equation} \label{eq:inhom_bil}
            \vn{ \Delta e^{r\Delta} u \, e^{s\Delta} v}_{L^2_{t,x}} \ls r^{-\frac14} s^{-\frac{n-1}{4}} \vn{u}_{X^{1,b}} \vn{v}_{X^{0,b}}.
    \end{equation}   
Additionally, the estimates above hold when $ u, v $ on the left are replaced by $ \bar{u}, v $ or $ u, \bar{v} $. 
\end{corollary}

The products above express high-low interactions with frequency separation as follows. 
If $ r = N^{-2},\ s = M^{-2} $, then heuristically $ e^{s\Delta} v \approx P_{<M}v $ and $ (1-e^{r\Delta})u \approx P_{>N} u $, 
while $ r \Delta e^{r\Delta} u \approx P_N u$. 
%Note that $ r^{\frac14} s^{-\frac{n-1}{4}} = N^{-\frac12} M^{\frac{n-1}{2}} $.
Moreover, it can be shown that the $ (-r \Delta)^{\theta} $ in \eqref{eq:heat_bil} can be replaced by $ r^{\frac12}\nabla$.

\subsection{Application: the cubic defocusing NLS below energy on $ \mathbb{H}^2 $} \

Consider the cubic NLS \eqref{eq:nls} on $ \mathcal{M}=\mathbb{H}^2 $ with initial data $ u(0,x) = u_0(x) \in H^{\sigma}(\mathbb{H}^2) $.

% \begin{equation} 
% i\partial_t u + \Delta u = |u|^2 u \quad \text{on } \mathbb{H}^2
% \end{equation}

Solutions of \eqref{eq:nls} enjoy conservation of mass $ \vn{u(t,\cdot)}_{L^2(\mathbb{H}^2)} = \vn{u_0}_{L^2(\mathbb{H}^2)} $ and energy (if finite)
\begin{equation*}
E(u)(t) = \int_{\mathbb{H}^2} \frac{1}{2} |\nabla u(t,x)|^2 + \frac{1}{4} |u(t,x)|^4 \, dx = E(u)(0).
\end{equation*}

At the energy regularity $ u_0 \in H^1(\mathbb{H}^2) $, Anker-Pierfelice \cite{anker2009nonlinear} and Ionescu-Staffilani \cite{ionescu2009semilinear} obtain global well-posedness as a consequence of energy conservation and the local well-posedness that they establish by proving Strichartz estimates.
The results in \cite{anker2009nonlinear, ionescu2009semilinear} apply to general power-type nonlinearity and also include scattering in $ H^1$. 

Below the energy regularity, for radial data in $ H^{\sigma}(\mathbb{H}^2) $, Staffilani-Yu
\cite{staffilani2020high} prove global well-posedness and scattering for $ \sigma > \frac{3}{4} $ by implementing the high-low method \cite{bourgain1998refinements} on $ \mathbb{H}^2 $. See also \cite{ma2022scattering} for radial $ \mathbb{H}^3 $.

The goal here is to relax the radial assumption and the regularity condition of the initial data that ensures global solutions 
by adapting the classical $I$-method of \cite{colliander2002almost} to the hyperbolic setting.
% of Colliander-Keel-Staffilani-Takaoka-Tao .

\begin{theorem} \label{thm:GWP_NLS}
The initial value problem \eqref{eq:nls} is globally well-posed for initial data $ u_0 \in H^{\sigma}(\mathbb{H}^2) $ for $ \sigma > \frac{4}{7} $.
\end{theorem}
This means that for any $ T>0 $ there is a unique solution $ u \in C([0,T], H^{\sigma}(\mathbb{H}^2)) $ which depends continuously on $ u_0\in H^{\sigma}(\mathbb{H}^2) $.
The proof provides a polynomial bound $ \|u(T)\|_{H^\sigma} \lesssim_{\vn{u_0}_{H^{\sigma}}} (1+T)^{C_{\sigma}}$. 

\

Following the $I$-method philosophy, the main ingredient is proving an almost-conservation law (Prop. \ref{prop:almost_conservation}) for a modified energy functional, which here is of type
\begin{equation*}
 (N^2 s)^{1-\sigma} E(e^{s\Delta}u)(t), 
\end{equation*}
for $ s \le N^{-2}$. 
See Definition \ref{def:modified_energy} for a complete definition, remark \ref{rmk:modified_energy_heuristics} for motivation, and remark \ref{rmk:sim_diff_eucl} for a comparison with the euclidean \cite{colliander2002almost}.
This heat-based approach to the $I$-method was introduced by the author in \cite{gavrus2022gauge} in the context of formulating a gauge-invariant modified energy for the Yang-Mills equation. 

\

To put Theorem \ref{thm:GWP_NLS} into context note that the condition $ \sigma > \frac{4}{7} $ mirrors the regularity threshold of the “first generation" $I$-method paper \cite{colliander2002almost}.   
Subsequent works on $ \mb{R}^2 $ lowered the regularity using refinements such as: 
multilinear correction (\cite{colliander2008resonant}, $ \sigma > \frac{1}{2} $), 
combination with interaction Morawetz (\cite{colliander2007interactionMorawetz}, $ \sigma > \frac{2}{5} $), 
bootstrapped Morawetz (\cite{colliander2010bootstrapped}, $ \sigma > \frac{1}{3} $), 
improved almost Morawetz (\cite{dodson2011improved}, $ \sigma > \frac{1}{4} $). 
The $I$-method loses strength as $ \sigma$ approaches $ 0 $.
The $ L^2_x(\mathbb{R}^2) $ endpoint was reached by \cite{killip2009cubic} for radial data, culminating with Dodson \cite{dodson2016global} removing the radial assumption.

We implemented the hyperbolic space $I$-method on the model of the early \cite{colliander2002almost} because
it provides a clear and concise nonlinear application of the bilinear Strichartz estimate in Thm \ref{thm:bilinear-hyperbolic}, beyond what linear norms alone can achieve
and it relaxes the regularity requirements of existing results \cite{anker2009nonlinear,ionescu2009semilinear,staffilani2020high}.
 
Moreover, it provides a second example of bypassing the lost euclidean relation between physical space multiplication and convolution in Fourier space, which becomes a unifying theme of this paper.
Instead, it is shown that a modified energy and its commutator structure can be addressed entirely in physical space. 
The Fourier multiplier commutator of the $\mb{R}^2$ modified energy is replaced by the heat–Schr\"odinger commutator $ |u_s|^2 u_s - e^{s\Delta}(|u|^2 u) $, where $u_s=e^{s\Delta}u$, whose 
defect formula is expressed through terms like $ \langle \nabla u_r, \nabla u_r \rangle \bar{u}_r $
realizing the cancellation geometrically rather than symbolically, see remark \ref{rmk:sim_diff_eucl} for a detailed comparison. 
After combining this with the bilinear Strichartz estimate \eqref{eq:inhom_bil}, 
one recovers the same effective quartic numerology as in \cite{colliander2002almost}'s $\mb{R}^2$ argument.

\

\noindent {\bf Related works.}
Dispersive and Strichartz estimates on $ \Hn$ are established in \cite{anker2009nonlinear,banica2007nonlinear,ionescu2009semilinear}.
Papers on NLS on $ \Hn$ include \cite{banica2007nonlinear,banica2008scattering,banica2009scattering,ionescu2009semilinear,anker2009nonlinear,christianson2010existence,ionescu2012global,banica2015global,staffilani2020high,ma2022scattering,Wilson2023GlobalWA}.
Local smoothing on $ \Hn$ is studied in \cite{lawrie2023localsmoothing}.
Results on geometric equations on $ \Hn$ include \cite{lawrie2018cauchy,li2020asymptotic,li2022global,lawrie2023asymptotic}.
We mention also \cite{tataru2001strichartz,metcalfe2012dispersive,kaizuka2014resolvent,germain2023spectral}.
Results on NLS on compact manifolds include \cite{burq2005bilinear,herr2011global,ionescu2012energy,herr2013quintic,herr2026global}.
For the $I$-method on compact manifolds see \cite{hani2012global,de2007global}.

%\cite{tataru2001strichartz}
%\cite{banica2007nonlinear}
%\cite{banica2008scattering}
%\cite{anker2009nonlinear} 
%\cite{ionescu2009semilinear}
% \cite{metcalfe2012dispersive}
%\cite{ionescu2012global}
%\cite{kaizuka2014resolvent}
%\cite{banica2015global}
%\cite{staffilani2020high}
%\cite{lawrie2018cauchy}
%\cite{li2020asymptotic}
%\cite{li2022global}
%\cite{lawrie2023localsmoothing}
%\cite{lawrie2023asymptotic}
%\cite{Wilson2023GlobalWA}
%\cite{ma2022scattering}
%\qquad
% \cite{herr2011global}
% \cite{ionescu2012energy}
% \cite{herr2026global}
% \cite{herr2013quintic} \ 
% \cite{deng2025sharp}
%\qquad 
%\cite{hani2012global}
%\cite{de2007global}

\subsection{AI statement}  
This paper developed through collaboration between the author and AI (mostly GPT-5.x) in the first half of 2026. 

\begin{enumerate}
\item On Theorem \ref{thm:bilinear-hyperbolic}: 
the initial approach 
revisited past attempts of the author (from his postdoctoral period) to prove  \eqref{BS}
by extending the euclidean ``lateral estimates'' (see \cite{ionescu2006low,tao2010physical}) to $\mathbb{H}^n$. 
%(i.e. splitting \eqref{BS} into a lateral smoothing and a maximal estimate).
During a series of technical and exploratory conversations with GPT, 
the simpler product-space $TT^*$ approach emerged. 
In particular, a GPT search provided the identities from \cite{Koornwinder2018DualAddition,HallnasRuijsenaars2018Product}
that form the basis of Lemma \ref{lemma:bound:T}.
\item On Theorem \ref{thm:GWP_NLS}: the idea of applying the bilinear estimate \eqref{BS} to the cubic NLS on $\mathbb{H}^2$
is due to the author (following \cite{colliander2002almost} on $\mathbb{R}^2$).
The formulation of the modified energy in Def. \ref{def:modified_energy} 
comes from the author's \cite{gavrus2022gauge}, while the structuring of the argument follows the pattern of \cite{colliander2002almost}.
GPT helped with the numerology. 
\end{enumerate}

The article was written by the author, who takes full responsibility for its content and correctness.

%% file: sections/proof_bilinear.tex
\addtocontents{toc}{\protect\setcounter{tocdepth}{1}}
\section{Proof of Theorem \ref{thm:bilinear-hyperbolic} and Corollary \ref{cor:bilinear_strichartz_variations}}

%%%%%%%%%%%%%%%%%%%%%%%%%%%%%%%%%%%%%%%%%%%%%%%%%%%%%%%%%%%%%%%%%%%%%%%%%%%%%%
%%%%%%%%%%%%%%%%%%%%%%%%%%%%%%%%%%%%%%%%%%%%%%%%%%%%%%%%%%%%%%%%%%%%%%%%%%%%%%
%%%%%%%%%%%%%%%%%%%%%%%%%%%%%%%%%%%%%%%%%%%%%%%%%%%%%%%%%%%%%%%%%%%%%%%%%%%%%%

\subsection{Proof of Theorem \ref{thm:bilinear-hyperbolic}}

In the case $ N \simeq M $ the estimate follows from the Strichartz estimate $ e^{it\Delta}: L^2_x(\mathbb{H}^n) \to L^4_t L^q_x(\mathbb{R} \times \mathbb{H}^n) $ from \cite{anker2009nonlinear,ionescu2009semilinear}, for $ q = \frac{2n}{n-1} $ and the Sobolev embedding $ W^{\frac{n}{2}-1,q} \subset L^{2n}$. 

It remains to consider the separated frequencies case $ N \gg M $. 
Conjugating if needed, it suffices to consider the cases $ uv $ and $ u \bar{v} $.
Due to the spectral support assumption, one may apply appropriate projections $ P_N, P_M $ given by enlarged dyadic cutoffs multipliers
such that $ f=P_N f $, $ g=P_M g $.

Let $K_M $, $K_M^{(2)} $ be the kernels of $e^{\pm it\Delta}P_M$ and $e^{\pm it\Delta}P_M^2$, see \eqref{eq:convolution_def}, \eqref{eq:spherical_transform}, \eqref{eq:transformed_convolution} meaning
\begin{equation}\label{eq:kernel_pn}
K_M^{(2)}(t,r) = 
 \int_{0}^{\infty}
e^{\mp it(\mu^2+\rho^2)}\,
\chi_M^2(\mu)
\varphi_\mu(r)\,
|c(\mu)|^{-2}\,
d\mu,
\qquad r \ge 0. 
\end{equation}
where the $ \pm=+ $ for $ v $ and $ \pm=- $ for $ \bar{v} $. For $ K_N $ the sign is $ + $, since $ u $ is not conjugated.
%i.e.
%\[ \bigl(e^{it\Delta}P_N f\bigr)(x) = 
%\int_{\mathbb{H}^n} K_N(t;x,y)\,f(y)\,dy. \]

Consider the linear operator $ F(x,y)\mapsto e^{it\Delta_x}e^{\pm it\Delta_y}  F(x,x) $ on $L^2(\mathbb{H}^n\times\mathbb{H}^n)$, i.e. 
define 
\[
(TF)(t,x)
:=
\iint_{\mathbb{H}^n\times\mathbb{H}^n}
K_N(t, d(x,y))\,K_M(t, d(x,z))\,F(y,z)\,dy\,dz.
\]

For $F=f\otimes g$ one has $ T(f\otimes g)=uv $, otherwise $ T(f\otimes \bar{g})=u\bar{v} $ if the sign is $ \pm=-$. 

Then \eqref{BS} follows from the boundedness of the operator 
$$ M^{-\frac{n-1}{2}}\,N^{\frac12} T : L^2(\mathbb{H}^n\times\mathbb{H}^n)\to L^2(\mathbb{R}\times\mathbb{H}^n),$$
for which bounding the norm of $ T T^*: L^2(\mathbb{R}\times\mathbb{H}^n) \to L^2(\mathbb{R}\times\mathbb{H}^n)$ suffices. 
The kernel of $ T T^*$ is
\[
K_N^{(2)}(t-s, d(x,y)) \cdot K_M^{(2)}(t-s, d(x,y)) 
%\qquad \qquad (t,x), (s,y) \in \mathbb{R}\times\mathbb{H}^n.
\]
Therefore, $ T T^*$ is diagonalizable by Fourier transform in $t$ and spherical (radial) spectral decomposition in $ x $, see \eqref{eq:FT}, \eqref{eq:space_time_multiplier}. 
Writing
\[
\mathcal{F}_t \widetilde{T T^* F}(\tau,\nu, b) = m(\tau, \nu) \mathcal{F}_t \widetilde{F}(\tau,\nu, b)
\]
it suffices to bound $ m(\tau, \nu) $ uniformly. 

Using polar coordinates $dx = (\sinh r)^{n-1}\,dr\,d\theta$ and \eqref{eq:kernel_pn}, up to constants
\begin{align*}
 \mathcal{F}_t^{-1} m(t,\nu) & = \int_0^{\infty} K_N^{(2)}(t,r) K_M^{(2)}(t,r) \varphi_\nu(r)\,(\sinh r)^{n-1}\,dr \\
& = \iint e^{-it(\lambda^2 \pm \mu^2+\rho^2\pm \rho^2)}\,
\chi_N^2(\lambda) \chi_M^2(\mu)  |c(\lambda)|^{-2}\, |c(\mu)|^{-2} \mathcal{T}(\lambda,\mu,\nu)
d\lambda \,d\mu
\end{align*}
where $ \mathcal{T} $ is the Jacobi triple product (linearization coefficient) 
\begin{equation}\label{def:T:coeff}
\mathcal{T}(\lambda,\mu,\nu) := \widetilde{\varphi_\lambda\, \varphi_\mu}(\nu) = \int_0^{\infty} \varphi_\lambda(r) \varphi_\mu(r) \varphi_\nu(r)\,(\sinh r)^{n-1}\,dr 
%= b(\lambda,\mu,\nu) |c(\nu)|^{2}
\end{equation}
Therefore 
\begin{equation}\label{eq:m:integral}
m(\tau, \nu) = \iint \delta(\tau+\lambda^2 \pm \mu^2+C_{\rho})\, \chi_N^2(\lambda)\, \chi_M^2(\mu)\, \mathcal{T}(\lambda,\mu,\nu) \frac{d\lambda}{|c(\lambda)|^{2}} \frac{d\mu}{|c(\mu)|^{2}},
\end{equation}
where $ C_{\rho}= \rho^2\pm \rho^2$ is a harmless spectral-gap constant.
\begin{lemma}\label{lemma:bound:T} 
Denoting $ {\bf d}=\text{dist}\bigl(\nu,\,[|\lambda-\mu|,\lambda+\mu]\bigr) $ one has
\begin{equation}\label{eq:T-bound}
\mathcal{T}(\lambda,\mu,\nu) 
%= |c(\nu)|^{2} b(\lambda,\mu,\nu) 
\simeq
\frac{\,\displaystyle \prod_{\pm_1, \pm_2} \bigl(1+|\lambda\pm_1\mu\pm_2\nu|\bigr)^{\rho-1}}
{(1+\lambda)^{2\rho-1}(1+\mu)^{2\rho-1}(1+\nu)^{2\rho-1}}\,
e^{-\pi\, {\bf d}}.
%e^{-\pi\,\text{dist}\bigl(\nu,\,[|\lambda-\mu|,\lambda+\mu]\bigr)}.
\end{equation}
\end{lemma}
When $ n=3 $ \eqref{eq:T-bound} is elementary, see Remark \ref{remark:3d:Tbound}. On the support of the cutoffs $ \lambda \sim N, \ \mu \sim M $ and 
$ \mathcal{T}(\lambda,\mu,\nu) $ is exponentially small unless $ \nu \sim N $ and $ \nu \in [\lambda-\mu, \lambda+\mu] $ up to $ O(1)$. 

%%%%%%%%%%%%%%%%%%%%%%%%%%%%%%%%%%%%%%%%%%%%%%%%%%%%%%%%%%%%%%%%%%%%%%%%%%%%%%
%%%%%%%%%%%%%%%%%%%%%%%%%%%%%%%%%%%%%%%%%%%%%%%%%%%%%%%%%%%%%%%%%%%%%%%%%%%%%%
%%%%%%%%%%%%%%%%%%%%%%%%%%%%%%%%%%%%%%%%%%%%%%%%%%%%%%%%%%%%%%%%%%%%%%%%%%%%%%

{\bf (1) Case $ n \geq 3$.} Integrate in $ \lambda $ first and use the identity %$ 2 \delta \left(x^{2}-\alpha ^{2}\right)=|\alpha|^{-1} \big[\delta \left(x+\alpha \right)+\delta \left(x-\alpha \right) \big] $ 
$ 2 \delta \left(x^{2}-a \right)=a^{-1/2} \delta \left(x- \sqrt{a} \right)  $ on $ x>0 $, 
which provides a power $ N^{-1}$. 
Integration in $ \mu $ provides a power of $ M $. Use \eqref{eq:c:function} for the two $ |c|^{-2} $ and \eqref{eq:T-bound} for $ \mathcal{T} $. 
When $ \nu \in [\lambda-\mu, \lambda+\mu] $ two of the numerator terms $ |\lambda\pm_1\mu\pm_2\nu| $  are $ \sim N $ and two are $ O(M) $,
so $ \mathcal{T} \lesssim M^{-1} N^{-(n-1)} $. When $ \nu \notin [\lambda-\mu, \lambda+\mu]$ use the exponential decay as needed. 
When a power of $ N $ is missing, one recovers it either from $ (1+\nu)^{1-2\rho} $ or from $ e^{-\pi d} $. In all cases one obtains
\begin{equation}\label{eq:m:bound}
\|T\|^2 = \|TT^*\|_{L^2_{t,x}\to L^2_{t,x}} = \sup_{\tau,\nu} | m(\tau, \nu)| \lesssim M^{n-1} N^{-1}
\end{equation}
  
{\bf (2) Case $ n = 2$.} The exponents flip sign. The focus is on $ \nu \in [\lambda-\mu, \lambda+\mu] $. Then \eqref{eq:T-bound} becomes
\[
\mathcal{T}(\lambda,\mu,\nu) \ls N^{-1} \bigl(1+|\lambda+\mu-\nu|\bigr)^{-\frac{1}{2}} \bigl(1+|\lambda-\mu-\nu|\bigr)^{-\frac{1}{2}}.
\]   
Bounding this uniformly would lose a power $ M^{-\frac{1}{2}} $ near the endpoints $\nu = \lambda \pm \mu $, so instead this term must be integrated.  
The main observation is that in both $\pm $ cases the delta curve (circle or hyperbola) remains transverse to the endpoint singularity lines,
since for $ N \gg M $ the slope is almost vertical in the $ (\lambda,\mu) $ plane.
Split \eqref{eq:m:integral} according to the regions: $ \lmd < \nu $ and $ \lmd > \nu $. 
Without loss of generality consider $ \lmd < \nu $, in which case $ |\lambda-\mu-\nu|>M $.
Apply a change of variables: one may use polar coordinates, or set $ w=\lambda+\mu-\nu$. By the coarea formula 
$ \int f(x) \delta(g(x)) \, dx = \int_{g^{-1}(0)} f(x)/|\nabla g| \, d\sigma(x)$, 
the integral becomes over a curve of length $ O(M) $ and things reduce to bounding 
$ \int_{|w|\ls M} (1+|w|)^{-\frac{1}{2}} dw \ls M^{\frac{1}{2}} $, concluding \eqref{eq:m:bound}.
$ \qed $

%%%%%%%%%%%%%%%%%%%%%%%%%%%%%%%%%%%%%%%%%%%%%%%%%%%%%%%%%%%%%%%%%%%%%%%%%%%%%%
%%%%%%%%%%%%%%%%%%%%%%%%%%%%%%%%%%%%%%%%%%%%%%%%%%%%%%%%%%%%%%%%%%%%%%%%%%%%%%
%%%%%%%%%%%%%%%%%%%%%%%%%%%%%%%%%%%%%%%%%%%%%%%%%%%%%%%%%%%%%%%%%%%%%%%%%%%%%%

\subsection{Proof of Lemma \ref{lemma:bound:T}} \

To complete the proof of \eqref{eq:m:bound} %Theorem \ref{thm:bilinear-hyperbolic}, 
it remains to 
establish the bound \eqref{eq:T-bound}. %on $ \mathcal{T}(\lambda,\mu,\nu) $. 
Writing 
\begin{equation*} %\label{def:b}
\varphi_\lambda(r) \varphi_\mu(r) = \int_0^{\infty}  \varphi_\nu(r)\, b(\lambda,\mu,\nu) d \nu
\end{equation*}
then for a $ C_\rho > 0 $ one has 
\begin{equation} \label{eq:b:explicit}
b(\lambda,\mu,\nu)= C_\rho\,
\frac{\displaystyle \prod_{\pm_1, \pm_2}
\left|\Gamma\!\left(\frac{\rho+i(\lambda\pm_1\mu\pm_2\nu)}{2}\right)\right|^{2}}
{\displaystyle |\Gamma(i\nu)|^{2}\,|\Gamma(\rho+i\lambda)|^{2}\,|\Gamma(\rho+i\mu)|^{2}}
\end{equation}
by following the identities from \cite{Koornwinder2018DualAddition} and \cite{HallnasRuijsenaars2018Product}.
Begin with \eqref{eq:Jacobi:Gauss} and \cite[(2.22)]{Koornwinder2018DualAddition}. 
Plug this into 
the dual product formula \cite[(137)]{HallnasRuijsenaars2018Product}, \cite[(6.1)]{Koornwinder2018DualAddition}, 
cancel and regroup some $ \Gamma$ factors and conclude \eqref{eq:b:explicit}. 

Noting that $ \mathcal{T}(\lambda,\mu,\nu) = |c(\nu)|^{2} b(\lambda,\mu,\nu) $ from \eqref{eq:spherical_transform},
the bound \eqref{eq:T-bound} follows from \eqref{eq:b:explicit}, by combining \eqref{eq:c:function}, the strip bound \cite[{(5.11.9)}]{NIST:DLMF}
\begin{equation*} %\label{eq:gamma-strip-bounds}
c_a(1+|y|)^{a-\frac12}e^{-\frac{\pi}{2}|y|}
\le |\Gamma(a+iy)|
\le C_a(1+|y|)^{a-\frac12}e^{-\frac{\pi}{2}|y|}, \qquad a>0, \ y\in\mathbb{R}.
\end{equation*}
and the identity
\begin{equation*} 
|\lmd+\mu+\nu|+|\lmd+\mu-\nu|+|\lmd-\mu+\nu|+|\lmd-\mu-\nu| = 2(\lmd+\mu+\nu)+2{\bf d}. \qquad \qed
\end{equation*}

%%%%%%%%%%%%%%%%%%%%%%%%%%%%%%%%%%%%%%%%%%%%%%%%%%%%%%%%%%%%%%%%%%%%%%%%%%%%%%
%%%%%%%%%%%%%%%%%%%%%%%%%%%%%%%%%%%%%%%%%%%%%%%%%%%%%%%%%%%%%%%%%%%%%%%%%%%%%%
%%%%%%%%%%%%%%%%%%%%%%%%%%%%%%%%%%%%%%%%%%%%%%%%%%%%%%%%%%%%%%%%%%%%%%%%%%%%%%

\begin{remark}\label{remark:3d:Tbound}
The case $ n = 3 $ is more explicit and elementary. By \eqref{eq:3d:spherical} one has 
\begin{equation}\label{eq:T-integral}
% b(\lambda,\mu,\nu) |c(\nu)|^{2} =  
 \mathcal{T}(\lambda,\mu,\nu) = \frac{1}{\lambda\mu\nu}
\int_{0}^{\infty}\frac{\sin(\lambda r)\sin(\mu r)\sin(\nu r)}{\sinh r}\,dr 
= \frac{\pi}{8 \lambda\mu\nu} J(\lambda,\mu,\nu)
\end{equation}
where, 
\begin{equation}\label{eq:tanh:repr}
J(\lambda,\mu,\nu) = 
\textstyle
\tanh\!\left(\frac{\pi(\nu+\lambda-\mu)}{2}\right)
+\tanh\!\left(\frac{\pi(\nu-\lambda+\mu)}{2}\right)
-\tanh\!\left(\frac{\pi(\nu+\lambda+\mu)}{2}\right)
-\tanh\!\left(\frac{\pi(\nu-\lambda-\mu)}{2}\right)
\end{equation}
based on the trigonometric identity
\begin{equation*}
4 \sin(\lambda r)\sin(\mu r)\sin(\nu r)
=
\sin\bigl((\nu+\lambda-\mu)r\bigr)
+\sin\bigl((\nu-\lambda+\mu)r\bigr)
-\sin\bigl((\nu+\lambda+\mu)r\bigr)
-\sin\bigl((\nu-\lambda-\mu)r\bigr)
\end{equation*}
and the classical integral 
$  \int_{0}^{\infty}\frac{\sin(\sigma r)}{\sinh r}\,dr = \frac{\pi}{2}\tanh\!\left(\frac{\pi\sigma}{2}\right) $.
By the inequalities 
\[
1-\tanh x \le 2e^{-2x},
\qquad
1+\tanh(-x)\le 2e^{-2x}, \qquad x \geq 0.
\]
each $ \tanh $ term in \eqref{eq:tanh:repr} differs from $ \text{sign}(\eta)$ by $ O(e^{-c \eta})$ where $ \eta=\nu\pm\lambda\pm\mu $ and  one obtains 
\begin{equation*} %\label{eq:J-exp-decay}
|J(\lambda,\mu,\nu)|
\lesssim
e^{-c\,{\bf d} }, \qquad {\bf d} = \text{dist}\!\bigl(\nu,\,[|\lambda-\mu|,\lambda+\mu]\bigr).
\end{equation*}
Also $ J(\lambda,\mu,\nu) = O(\nu e^{-c\,{\bf d} })$ as $ \nu \to 0$, so \eqref{eq:T-integral} doesn't blow up. 
\end{remark}

%%%%%%%%%%%%%%%%%%%%%%%%%%%%%%%%%%%%%%%%%%%%%%%%%%%%%%%%%%%%%%%%%%%%%%%%%%%%%%
%%%%%%%%%%%%%%%%%%%%%%%%%%%%%%%%%%%%%%%%%%%%%%%%%%%%%%%%%%%%%%%%%%%%%%%%%%%%%%
%%%%%%%%%%%%%%%%%%%%%%%%%%%%%%%%%%%%%%%%%%%%%%%%%%%%%%%%%%%%%%%%%%%%%%%%%%%%%%

\subsection{Proof of sharpness} 
To complete the proof of Theorem \ref{thm:bilinear-hyperbolic}, it remains to prove that the bound in \eqref{BS} is optimal.
% The bound in \eqref{eq:m:bound} is actually $ \simeq M^{n-1} N^{-1} $ and one has 
%     $$ \sup_{\vn{f}=\vn{g}=1} \vn{T(f\otimes g)}_{L^2_{t,x}} \leq \vn{T}_{L^2_{x,y}\to L^2_{t,x}} \simeq M^{n-1} N^{-1}. $$
% The reverse inequality is false in general, particularly when the maximimizing tensor is entangled, not close to a product tensor.
For $ 1 \ls M \leq N $, let $A = \frac{M^2}{N}$. Choose central frequencies $ \lambda_0=\nu_0 \sim N, \ \mu_0 \sim M $, $ \tau_0 = -\lambda_0^2 \mp \mu_0^2-C_{\rho} $ and define intervals
\[ I_N = [\lambda_0 - cA, \lambda_0 + cA], \ I_M = [\mu_0 - cM, \mu_0 + cM], \ J_{\nu}= [\nu_0 - cM, \nu_0 + cM], \ J_{\tau}=[\tau_0 - cM^2, \tau_0 + cM^2]
\]
where $ c > 0 $. Let $ f, g $ be $ L^2 $-normalized radial initial data with spherical transforms
\[ \tilde{f}(\lambda) = A^{-\frac12} \vm{c(\lambda)} 1_{I_N}(\lambda), \qquad \tilde{g}(\mu) = M^{-\frac12} \vm{c(\mu)}  1_{I_M}(\mu).
\]
Then 
\begin{equation*}
\mathcal{F}_t \widetilde{uv}(\tau,\nu) = \int \delta(\tau+\lambda^2 \pm \mu^2+C_{\rho})\,  \tilde{f}(\lambda)\,  \tilde{g}(\mu)\,  \mathcal{T}(\lambda,\mu,\nu)\,  \frac{d\lambda}{|c(\lambda)|^2} \frac{d\mu}{|c(\mu)|^2}
\end{equation*}
On the chosen region the integrand is nonnegative and has fixed size, as follows.
On the region $ \lmd \in I_N, \ \mu \in I_M, \ \nu \in J_{\nu} $, one has $ \nu \in [\lambda-\mu, \lambda+\mu] $ up to harmless fixed constants, 
so the exponential factor in $\mathcal{T}(\lambda,\mu,\nu) $ is not small and $ \mathcal{T}(\lambda,\mu,\nu) \simeq M^{-1} N^{-(n-1)} $. 
For $(\tau,\nu) \in J_{\tau} \times J_{\nu}  $ the delta curve intersects $ I_N \times I_M $ in a set of $\mu$-length $ \simeq M $, by shrinking the constants if needed.
By integration it follows that $ |\mathcal{F}_t \widetilde{uv}(\tau,\nu)|^2 \gtrsim N^{-n} M^{n-4} $ on $J_{\tau} \times J_{\nu} $ and therefore $ \|uv\|_{L^2_{t,x}}^2 \gtrsim N^{-1} M^{n-1} $. \qed

\begin{remark}
The proof of \eqref{BS} provides additional information.
By applying an output projection $P_L$ with $L$ dyadically separated from $N$ 
%(say $L \le N/8$ or $L \ge 8N$), 
then an exponential gain absorbs the other powers:
\begin{equation*}
\|P_L(uv)\|_{L^2_{t,x}} \lesssim e^{-c \max(L,N)} \|f\|_{L^2_x} \|g\|_{L^2_x}, \qquad  
M\ll N; L \ll N \ \text{or} \ L \gg N. 
\end{equation*}
This relates to an approximate version of the euclidean Littlewood-Paley trichotomy. 
\end{remark}

%%%%%%%%%%%%%%%%%%%%%%%%%%%%%%%%%%%%%%%%%%%%%%%%%%%%%%%%%%%%%%%%%%%%%%%%%%%%%%
%%%%%%%%%%%%%%%%%%%%%%%%%%%%%%%%%%%%%%%%%%%%%%%%%%%%%%%%%%%%%%%%%%%%%%%%%%%%%%
%%%%%%%%%%%%%%%%%%%%%%%%%%%%%%%%%%%%%%%%%%%%%%%%%%%%%%%%%%%%%%%%%%%%%%%%%%%%%%

\subsection{Proof of Corollary \ref{cor:bilinear_strichartz_variations}} 

{\bf(1)} Normalize $\|f\|_{L^2} =\|g\|_{L^2}=1$. Split
\begin{equation*}
    (-r \Delta)^{\theta} e^{r\Delta} u = \sum_{N\geq 1} a_N u_N, \qquad  e^{s\Delta} v = \sum_{M \geq 1} b_M v_M; \qquad a_N = (r N^2)^{\theta} e^{-cr N^2}, \ b_M = e^{-cs M^2}
\end{equation*}
where $ u_N = e^{it\Delta} f_N $, $ v_M = e^{it\Delta} g_M $ with $ f_N, g_M $ spectrally supported at dyadic frequencies $ N, M $ respectively and $ \vn{f_N}_{L^2}, \vn{g_M}_{L^2} \ls 1 $.
Then, by Theorem \ref{thm:bilinear-hyperbolic} 
\begin{equation*}
\vn{(-r \Delta)^{\theta} e^{r\Delta} u \, e^{s\Delta} v}_{L^2_{t,x}} 
\ls \sum_{N,M} a_N b_M M^{\frac{n-1}{2}}N^{-\frac12} \ls r^{\frac14} s^{-\frac{n-1}{4}} 
\end{equation*}
since 
\[
\frac{\min(M,N)^{\frac{n-1}{2}}}{\max(M,N)^{\frac12}} \leq M^{\frac{n-1}{2}}N^{-\frac12}, \qquad \sum_N N^{-\frac12} a_N \simeq r^{\frac14}, \qquad \sum_M M^{\frac{n-1}{2}} b_M \simeq s^{-\frac{n-1}{4}}.
\]
The proof for $ (1-e^{r\Delta})u $ is similar, using coefficients $ c_N \leq 1-e^{-cr N^2} \ls \min(1, r N^2) $ instead of $ a_N $.

{\bf(2)} The inhomogeneous estimate \eqref{eq:inhom_bil} follows from \eqref{eq:heat_bil} with $ \theta=\frac{1}{2} $ using the transfer principle for $ X^{s,b} $ spaces.
This means decomposing $ u, v $ using \eqref{eq:integral_rep} as
\[
u(t) = \int_{\mathbb{R}} e^{ith_1} e^{it\Delta} f_{h_1} \dd h_1, \qquad v(t) = \int_{\mathbb{R}} e^{ith_2} e^{it\Delta} f_{h_2} \dd h_2 
\]
applying \eqref{eq:heat_bil} to $ (-\Delta)^{\frac{1}{2}} f_{h_1}, f_{h_2} $ and concluding with Minkowski's inequality in $ h_1, h_2 $.

\

%%%%%%%%%%%%%%%%%%%%%%%%%%%%%%%%%%%%%%%%%%%%%%%%%%%%%%%%%%%%%%%%%%%%%%%%%%%%%%
%%%%%%%%%%%%%%%%%%%%%%%%%%%%%%%%%%%%%%%%%%%%%%%%%%%%%%%%%%%%%%%%%%%%%%%%%%%%%%
%%%%%%%%%%%%%%%%%%%%%%%%%%%%%%%%%%%%%%%%%%%%%%%%%%%%%%%%%%%%%%%%%%%%%%%%%%%%%%

%% file: sections/ac_main_thm_proof.tex
\section{Almost-conservation and Proof of Theorem \ref{thm:GWP_NLS}} 

In this section the modified energy and the heat-based I-method in physical space are introduced on $ \mathbb{H}^2 $,
leading to a reduction of Theorem \ref{thm:GWP_NLS} to an almost conservation law for the modified energy, stated in Prop. \ref{prop:almost_conservation}. 
This heat-based approach to the I-method was introduced by the author in \cite{gavrus2022gauge} in the context of formulating a gauge-invariant modified energy for the Yang-Mills equation. 

\begin{definition} \label{def:modified_energy}
Let $s>0$. Define the heat-smoothed field and the energy at heat scale $s$ by:
\begin{equation}
u_s(t,\cdot) := e^{s\Delta} u(t,\cdot), \qquad E(t,s) := E(u_s(t)) := \int_{\mathbb{H}^2} \frac{1}{2} |\nabla u_s(t,x)|^2 + \frac{1}{4} |u_s(t,x)|^4 \, dx.
\end{equation} 
Define the \emph{modified energy} by
\begin{equation}
\mathcal{IE}_N(t) := \sup_{0 < s \le s_0} (N^2 s)^{1-\sigma} E(t,s) + \int_0^{s_0} (N^2 s)^{1-\sigma} E(t,s)\,\frac{ds}{s}, \quad s_0 = N^{-2}.
\end{equation}
\end{definition}

\begin{remark}[Modified energy heuristics] \label{rmk:modified_energy_heuristics}
The classical euclidean I-method \cite{colliander2002almost} is based on $E[I_N\phi]$ involving the Fourier multiplier $ I = I_N $ with symbol
\begin{equation} \label{eq:I_symbol}
m(\xi) = 1 \quad \text{if} \quad |\xi| \leq N \qquad \text{and} \quad  (N/|\xi|)^{1-\sigma} \quad \text{if} \quad  |\xi| >2 N 
\end{equation} 
The Gaussian multipliers $ e^{s\Delta} $ concentrates the bulk of the frequencies of $ u_s $ to $ \jb{\lambda} \ls s^{-\frac12}$,
while $ \nabla e^{s\Delta} $ concentrates to $ \jb{\lambda} \simeq s^{-\frac12} $.
Denoting by $ x=\lambda^2+\frac14 $ the spectral variable of $ -\Delta_{\mathbb{H}^2} $, the gradient part of $ \mathcal{IE}_N $ has weight
\begin{equation} \label{eq:weight_def}
w(x)=\sup_{0 < s \le N^{-2}} (N^2 s)^{1-\sigma} x e^{-2sx} + \int_0^{N^{-2}} (N^2 s)^{1-\sigma} x e^{-2sx} \frac{ds}{s}.
\end{equation}
For $ x \ls N^2 $ ($ \lmd \ls N $), the weight is $ w(x) \simeq x $, while for $ x \gtrsim N^2 $, the weight is $ w(x) \simeq N^{2(1-\sigma)} x^\sigma $.
This matches the numerology \eqref{eq:I_symbol} (after taking the gradient).
The integral in \eqref{eq:weight_def} is the Littlewood–Paley square-function part.
\end{remark}

As in the euclidean I-method, the modified energy $ \mathcal{IE}_N(t) $ can be compared with Sobolev norms:

\begin{lemma} \label{lem:energy_comparison}
    For $ 0 < \sigma < 1 $, denoting $ P=e^{N^{-2} \Delta} $ one has
\begin{equation} \label{eq:energy_comparison_1}
    \vn{ u(t) }_{H^{\sigma}(\mathbb H^n)} \lesssim 
    \vn{P u(t)}_{H^1} + \vn{(1-P) u(t)}_{H^{\sigma}} \ls \mathcal{IE}_N(t)^{1/2}.
\end{equation}
\begin{equation} \label{eq:energy_comparison_2}
\mathcal{IE}_N(t)^{1/2} \ls N^{1-\sigma} \vn{u(t)}_{H^{\sigma}(\mathbb H^n)} + \vn{u(t)}_{L^4(\mathbb H^n)}^2.
\end{equation}
\end{lemma}
The proof is in section \ref{sec:proof_lemmas_prop}. The $ X^{s,b} $ spaces used below are defined in section \ref{sec:xsb_spaces} of the appendix.

\begin{definition} \label{def:Y_norm}
Given $ \sigma, b \in \mathbb{R} $, an interval $I $ and $ P = e^{N^{-2} \Delta} $, define the norm
 \begin{equation} \label{eq:Y_norm}
        \vn{u}_{Y(I)} := \vn{P u}_{X^{1,b}(I)} + N^{1-\sigma} \vn{(1-P) u}_{X^{\sigma,b}(I)}
\end{equation}
and, denoting $ u_s = e^{s\Delta} u $ and $ s_0 = N^{-2} $, the comparable norm
\begin{equation} \label{eq:tilde_Y_norm}
    \vn{u}_{\tilde{Y}(I)} := \|u_{s_0}\|_{X^{1,b}(I)} + \left( \int_0^{s_0} (N^2 s)^{1-\sigma} \| u_s\|_{X^{1,b}(I)}^2 \, \frac{ds}{s} \right)^{1/2} .
\end{equation}
\end{definition}
These norms will be used to control the solutions on small intervals. 

\begin{remark}
\begin{enumerate}
    \item The heat localization $ P = e^{N^{-2} \Delta} $ in \eqref{eq:Y_norm} can be replaced by a smooth spectral cutoff $ \chi(-\Delta/N^2) $
    and the resulting norm is equivalent.
    \item It is not hard to see that $ Y(I) = X^{1,b} + N^{\sigma-1} X^{\sigma,b} $ with 
    \begin{equation*}
     \vn{u}_{Y(I)} \simeq \inf_{u=u_1+u_2} \vn{u_1}_{X^{1,b}(I)} + N^{1-\sigma} \vn{u_2}_{X^{\sigma,b}(I)}.
    \end{equation*}
\end{enumerate}
\end{remark}

\begin{lemma} \label{lem:Y_equiv}
The norms in Definition \ref{def:Y_norm} are equivalent: $ \vn{u}_{\tilde{Y}(I)} \simeq \vn{u}_{Y(I)} $.
\end{lemma}
\begin{proof}
    When $ I= \mathbb{R}$ this follows from the symbol identity
\begin{equation} \label{eq:symbol_identity}
x e^{-2x/N^2} + \int_0^{N^{-2}} (N^2 s)^{1-\sigma} x e^{-2sx} \, \frac{ds}{s} \ \simeq \ x e^{-2x/N^2} + N^{2(1-\sigma)}x^\sigma(1 - e^{-x/N^2})^2,
\end{equation}
where $ x $ is the spectral variable of $ -\Delta_{\mathbb{H}^2} $. 
The general interval case follows by extension arguments. 
\end{proof}

In what follows the constants will be chosen such that 
\begin{equation} \label{eq:eta_def}
    0<\eta<\frac{7\sigma-4}{4(1-\sigma)}, \qquad \frac{1+\eta/4}{1+\eta} < 1- 2 \ep, \qquad \eta<1
\end{equation}

\begin{equation} \label{eq:NLnumerology}
N^{-2 (1-\sigma)} \ls L^{-1}, \qquad N^{\ep} \ls L^{\eta/8}
\end{equation}

\begin{proposition}[Modified LWP] \label{prop:modified_LWP}
 Let $ \sigma \in (\frac{4}{7}, 1) $, $ \eta, \ep > 0 $ satisfying \eqref{eq:eta_def},
 $ b = \frac{1}{2}+\ep $ and assume \eqref{eq:NLnumerology} holds. 
 Assume the inital data for \eqref{eq:nls} satisfies $ \mathcal{IE}_N(0) \le L $. 
 Then there exists $ \delta > 0 $ (depending on $ \sigma,\eta,\ep $ and the implicit constant in \eqref{eq:NLnumerology}) such that the maximal $ H^{\sigma} $ solution $ u $ of \eqref{eq:nls} exists on $ J = [0, \delta L^{-1-\eta}] $ 
 and satisfies $ \vn{u}_{Y(J)} \ls  L^{\frac12} $. \
 Moreover, one has $ \|u\|_{X^{0,b}(J)}\lesssim_{\|u_0\|_{L^2}} 1. $
\end{proposition}
The proof is relegated to section \ref{sec:proof_lemmas_prop}. This is the unscaled analogue of \cite[Prop 3.2]{colliander2002almost}, with equivalent numerology consistent with euclidean scaling, apart from small exponents.

\

\begin{proposition}[Almost conservation] \label{prop:almost_conservation}
    Let $ \sigma \in (\frac{4}{7}, 1) $ and assume \eqref{eq:NLnumerology} holds. 
    Assume the inital data for \eqref{eq:nls} at time $ t_0 $ is in $ C_c(\mathbb H^2) $ and  satisfies $ \mathcal{IE}_N(t_0) \le L $. 
    Then the solution $ u \in C(J, H^{\sigma}(\mathbb{H}^2))$ of \eqref{eq:nls} from Prop. \ref{prop:modified_LWP}, for $ J = [t_0, t_0+ \delta L^{-1-\eta}] $  satisfies
    \begin{equation} \label{eq:almost_conservation_s}
        (N^2 s)^{1-\sigma} \vm{E(t,s) - E(t_0,s)} \ls L \Big( \frac{L^{1/2}}{N} \Big)^{\frac{3}{2}} (N^2s)^{\frac{\sigma}{2}+\frac{1}{4}-} , 
     \end{equation}
     for all $ s \in (0, s_0] $ and $t \in J$, % for some fixed small $\epsilon>0$
    , where $s_0=N^{-2}$. Consequently,
    \begin{equation} \label{eq:almost_conservation}
        |\mathcal{IE}_N(t) - \mathcal{IE}_N(t_0)| \le C L \Big( \frac{L^{1/2}}{N} \Big)^{\frac{3}{2}}.
    \end{equation}
    where $ C $ depends on $ \vn{u(t_0)}_{L^2} $ and the implicit constant in \eqref{eq:NLnumerology}.
\end{proposition}
 The proof is deferred to section \ref{sec:proof_ac}. The power $\frac{3}{2} $ corresponds to the condition $ \sigma > \frac{4}{7} $. 

\begin{remark}[Similarity and differences with the $\mb{R}^2$ case in \cite{colliander2002almost}] \label{rmk:sim_diff_eucl} \ 
\begin{enumerate}
    \item The first step of the $ I$-method is usually to scale the modified energy to unit size, simplifying the numerology. 
Since scaling is not available on $ \mathbb{H}^n $, this step is skipped here, working directly with the original solution.
Still, the numerology mirrors \cite{colliander2002almost} with reversed scaling. 
    \item The smoothing in the $\mb{R}^2 $ energy $E[I_N\phi]$ uses a Fourier multiplier 
operator $ I=I_N $ with symbol \eqref{eq:I_symbol}. %$ m $ equal to $ 1 $ for frequencies $ \leq N $ and $ (N/|\xi| )^{1-\sigma} $ for higher frequencies. 
The hyperbolic smoothing is realized in a parabolic way. 
    \item The increment of the modified energy in \cite{colliander2002almost} contains a symbolic commutator $ [\mathcal{N}, I]$, where $ \mathcal{N} $ denotes taking the cubic nonlinearity. 
Here it becomes a heat-Schr\"odinger commutator $ [\mathcal{N}, e^{s\Delta}]$, measuring the failure of $e^{s\Delta}u$ to solve the cubic NLS \eqref{eq:nls}.
The main term in the \cite{colliander2002almost} energy increment is 
\begin{equation*}
\iint \overline{\Delta \phi} \big( \vm{I \phi}^2 I\phi - I(\vm{\phi}^2\phi) \big) dx\, dt
\end{equation*}
The euclidean cancelation is expressed in Fourier space through the symbol $ 1-\frac{m(\xi_2+\xi_3+\xi_4)}{m(\xi_2)m(\xi_3)m(\xi_4)}$. \
The hyperbolic commutator in the modified energy increment 
\begin{equation*}
\iint  \overline{\Delta u_s}\big(|u_s|^2 u_s - e^{s\Delta}(|u|^2 u) \big) dx\, dt,  % \qquad u_s = e^{s\Delta}u
\end{equation*}
is in physical space, containing terms $ r \langle \nabla u_r, \nabla u_r \rangle \bar{u}_r $, with $ u_r = e^{r\Delta}u$, $ r \leq s $, see Lemma \ref{lem:commutator_structure}. 
\item In \cite{colliander2002almost} the hard quartic case is the high-high-low-low interaction, 
where the mean value theorem for $ m $ provides a $ \frac{N_{low}}{N_{high}}$ gain. 
In the heat formulation, the low-frequency influence is suppressed by the gradient pair $ \langle \nabla u_r, \nabla u_r \rangle $.
This means that at most one term (the $ u_r $) can genuinely be low-frequency, 
in which case the bilinear Strichartz estimate for that part of the pair $ \Delta u_s u_r $ provides an additional gain due to frequency separation. 
\end{enumerate}
\end{remark}

\subsection{Reduction of Theorem \ref{thm:GWP_NLS} to the almost conservation law Prop. \ref{prop:almost_conservation}} \ 

Due to density, it can be assumed that $ u_0 $ is sufficiently regular. 
Let $u$ be the maximal $H^\sigma$ solution of \eqref{eq:nls} on its lifespan $[0,T_{\max})$, which  
exists due to the local well-posedness which follows from Strichartz estimates. 
To prove that $ T_{\max} $ is infinite
it suffices to prove a polynomial bound $ \vn{u(t)}_{H^{\sigma}(\mathbb{H}^2)} \le C_1 (1+t)^{K} $, where $ C_1 $ depends only on $ \vn{u_0}_{H^{\sigma}} $. 
The regularity assumption is removed by standard limiting arguments.

Fix $\sigma\in(\frac47,1)$ (for $\sigma\geq1$ see \cite{anker2009nonlinear}, \cite{ionescu2009semilinear}) 
and choose $ \eta \in (0,1) $ %, then choose $ \epsilon>0 $ sufficiently small 
such that 
\begin{equation}\label{eq:kappa_def}
0<\eta<\frac{7\sigma-4}{4(1-\sigma)}, \qquad \kappa:=\frac32-\Big(\frac72+2\eta\Big)(1-\sigma)>0
\end{equation}
 
By Lemma \ref{lem:energy_comparison} and Sobolev embedding $ H^{\sigma} \subset L^4 $ one has \ 
$ \mathcal{IE}_N(0)\lesssim_{\vn{u_0}_{H^{\sigma}}} N^{2(1-\sigma)}$. 
Choose $C_0=C_0(\vn{u_0}_{H^{\sigma}})\gg1$ and define
\begin{equation}\label{eq:LN_def}
L_N:=C_0\,N^{2(1-\sigma)} \qquad \text{so that }\qquad \mathcal{IE}_N(0)\le \frac12 L_N.
\end{equation}
For this choice the numerology \eqref{eq:NLnumerology} holds (after fixing a small $\varepsilon>0$ depending on $\sigma,\eta$).

Let $T>0$ be arbitrary. Apply Proposition \ref{prop:almost_conservation} and \ref{prop:modified_LWP} first with $ t_0 = 0 $ and then inductively on consecutive intervals of length at most $ \delta L_N^{-1-\eta} $
until reaching $ T $ while still maintaining the size of each $ \mathcal{IE}_N(t) \leq L_N $. 
Indeed, the number of intervals is $ O(T L_N^{1+\eta}) $ and each increment is bounded by \eqref{eq:almost_conservation} so the condition becomes 
\begin{equation*}
    C T L_N^{1+\eta} L_N \big(L_N^{1/2}/N \big)^{\frac{3}{2}} \leq L_N, \qquad \quad C=C(\vn{u_0}_{H^{\sigma}}), 
\end{equation*}
which, given the choices \eqref{eq:LN_def} and \eqref{eq:kappa_def}, becomes $ \tilde{C} T \leq N^{\kappa} $ for some $\tilde{C}=\tilde{C}(\vn{u_0}_{H^{\sigma}})$. 
This holds if $N$ is chosen sufficiently large depending on $T$ and $u_0$, with $ N \simeq (1+T)^{\frac{1}{\kappa}} $. 

Finally, Lemma \ref{lem:energy_comparison} and \eqref{eq:LN_def} ensure
\begin{equation*}
\|u(t)\|_{H^\sigma(\mathbb H^2)}
\lesssim \mathcal{IE}_N(t)^{1/2}
\lesssim L_N^{1/2} \lesssim_{\vn{u_0}_{H^{\sigma}}} N^{1-\sigma} \lesssim_{\vn{u_0}_{H^{\sigma}}} (1+T)^{\frac{1-\sigma}{\kappa}}. 
\end{equation*}
for $ t \in[0, T]$.   
Since $ T $ was arbitrary, it follows $ T_{\max}=\infty $,
concluding the proof of Theorem \ref{thm:GWP_NLS}.

\

%% file: sections/proof_ac_law.tex
\section{Proof of the almost conservation Proposition \ref{prop:almost_conservation}} \label{sec:proof_ac}

Recall Definition \ref{def:modified_energy} and differentiate 
\begin{equation} \label{eq:energy_derivative}
\frac{d}{dt}E(t,s) = \Re \int_{\mathbb{H}^2} \overline{\partial_t u_s} \left( -\Delta u_s + |u_s|^2 u_s \right) dx.
\end{equation}
Applying $e^{s\Delta}$ to \eqref{eq:nls} gives 
$ -\Delta u_s = i\partial_t u_s - e^{s\Delta}(|u|^2 u). $
Plugging this into \eqref{eq:energy_derivative} and integrating in time leads to
\begin{equation*}
E(t,s) - E(t_0,s) = \Re \int_{t_0}^{t} \int_{\mathbb{H}^2} \overline{\partial_t u_s}\, \mathcal{D}_s(u) \, dx \, dt.
\end{equation*}
where 
\begin{equation} \label{eq:commutator_def}
\mathcal{D}_s(u) := |u_s|^2 u_s - e^{s\Delta}(|u|^2 u).
\end{equation}
since the real part of $ i \vm{\partial_t u_s}^2 $ is $0$. Plugging the equation for $ \pt_t u_s $ splits the increment as  
\begin{equation*}
E(t,s) - E(t_0,s) = \mathrm{Term}_{1,s} + \mathrm{Term}_{2,s}
\end{equation*}
into a quartic and a sextic term 
\begin{equation} \label{eq:term1_term2_def}
\mathrm{Term}_{1,s} = \Im \int_{t_0}^{t} \int_{\mathbb{H}^2} \overline{\Delta u_s}\,\mathcal{D}_s(u)\,dx \, dt, \qquad 
\mathrm{Term}_{2,s} = - \Im \int_{t_0}^{t} \int_{\mathbb{H}^2} \overline{e^{s\Delta}(|u|^2 u)}\,\mathcal{D}_s(u)\,dx \, dt.
\end{equation}

The main claims are 
\begin{align} \label{eq:term1_bound}
|\mathrm{Term}_{1,s}| & \ls \frac{s^{\frac{3}{4}}}{(Ns^{\frac{1}{2}})^{3(1-\sigma)}} L^{\frac{7}{4}} \\ 
\label{eq:term2_bound}
|\mathrm{Term}_{2,s}| & \ls \frac{s}{(N^2s)^{2(1-\sigma)}} N^{0+} L^{2}
\end{align}
under the conditions of Prop. \ref{prop:modified_LWP}.
These imply \eqref{eq:almost_conservation_s} and,
after integration, together with $ N^2 s \leq 1 $, \eqref{eq:almost_conservation} follows.  

%The standard lemma below is a heat-flow version of the paraproduct decomposition of $\vm{u}^2 u$. 
\begin{lemma}[Commutator identity] \label{lem:commutator_structure}
%Denoting $ \Gamma(f,g) := \langle \nabla f, \nabla g \rangle_{\mathbb{H}^2} $, one has
The defect term $ \mathcal{D}_s(u) $ in \eqref{eq:commutator_def} is expressed as 
\begin{equation} \label{eq:commutator_structure}
\mathcal{D}_s(u) = -4 \int_0^s e^{(s-r)\Delta}( \langle \nabla u_r, \nabla \overline{u_r} \rangle u_r)\,dr
    - 2 \int_0^s e^{(s-r)\Delta}(\langle \nabla u_r, \nabla u_r \rangle \overline{u_r})\,dr.
\end{equation}
\end{lemma}
\begin{proof}
For general $f,g,h$ define $
G(r) := e^{(s-r)\Delta}\big((e^{r\Delta} f)(e^{r\Delta} g)(e^{r\Delta} h)\big),
$ so 
\begin{equation*}
(e^{s\Delta} f)(e^{s\Delta} g)(e^{s\Delta} h) - e^{s\Delta}(fgh) = \int_0^s G'(r)\, dr.
\end{equation*}
where $
G'(r) = e^{(s-r)\Delta}\big(-\Delta(uvw) + (\Delta u)vw + u(\Delta v)w + uv(\Delta w)\big),
$ with $u = e^{r\Delta} f$, $v = e^{r\Delta} g$, $w = e^{r\Delta} h$.
Expanding $ \Delta(uvw) $ using the Leibniz rule and collecting the terms gives \eqref{eq:commutator_structure}. 
\end{proof}

% \begin{remark}[Parallelism with the $\mb{R}^n$]
% By plugging in the representation \eqref{eq:commutator_structure} of $ \mathcal{D}_s(u) $ into \eqref{eq:term1_term2_def} one obtains a 
% \end{remark}

\subsection{Envelope bounds} 
The constants in this section may depend on $ \vn{u_0}_{L^2} $ but not on $ L, N $.
We are under the conditions of Prop. \ref{prop:modified_LWP}, so $ \vn{u}_{Y(J)} \ls L^{\frac12} $ and $ \vn{u}_{X^{0,b}} \ls 1$, so by \eqref{eq:strichartz_embedding},  
\begin{equation} \label{eq:envelope_bound_L4}
\vn{u_s}_{L^4_{t,x}} \ls \vn{u_s}_{X^{0,b}} \ls 1.
\end{equation}

Based on $ \partial_s u_s = \Delta u_s $, expanding
\begin{equation} \label{eq:heat_expansion}
u_s = u_{s_0} - \int_s^{s_0} r \Delta u_r \, \frac{dr}{r}, \qquad s \in (0, s_0],\ s_0=N^{-2}
\end{equation}
using Minkowski's inequality and Cauchy-Schwarz in $ dr/r $, from the definition \eqref{eq:tilde_Y_norm} it follows that
\begin{equation} \label{eq:envelope_bound_X1b}
\vn{u_s}_{X^{1,b}} \ls L^{\frac12}  (Ns^{\frac{1}{2}})^{-(1-\sigma)}, \qquad s \in (0, s_0].
\end{equation}
Since $ \vn{u_r}_{X^{1,b}} \ls \vn{u_{s_0}}_{X^{1,b}} \ls L^{\frac12} $, for $ r>s_0 $, the last two can be combined as
\begin{equation} \label{eq:envelope_bound_X1b_all_s}
\vn{u_s}_{X^{1,b}} \ls L^{\frac12} \max(1, (Ns^{\frac{1}{2}})^{-(1-\sigma)}), \qquad s > 0.
\end{equation}

As a consequence, using \eqref{eq:strichartz_embedding}
\begin{equation} \label{eq:env_nabla_ur_L4}
\vn{\nabla u_s}_{L^4_{t,x}} \ls  L^{\frac12} (Ns^{\frac{1}{2}})^{-(1-\sigma)},  \qquad s \in (0, s_0].
\end{equation}
 
By \eqref{eq:heat_expansion}, Sobolev $ H^{1+}\subset L^{\infty}_x $, Bernstein $ e^{\tau \Delta}: L^2_x \to \tau^{-\frac12} L^{\infty}_x $, and splitting $ u_r = e^{r/2 \Delta} u_{r/2} $
\begin{equation} \label{eq:envelope_bound_inf}
\vn{u_s}_{L^{\infty}_{t,x}} \ls N^{0+} L^{\frac12} (Ns^{1/2})^{-(1-\sigma)},  \qquad s \in (0, s_0].
\end{equation}

\subsection{Proof of the Term1 bound \eqref{eq:term1_bound}} \

By plugging the representation of $ \mathcal{D}_s(u) $ given by Lemma \ref{lem:commutator_structure} into $ \mathrm{Term}_{1,s} $, and using the self-adjointness of $ e^{(s-r)\Delta} $ 
it suffices to bound
\begin{equation} \label{eq:term1_bound_reduced}
    | \int_{t_0}^{t} \int_{\mathbb{H}^2}   e^{(2s-r)\Delta} \Delta \bar{u} \, \langle \nabla u_r \, ,   \nabla u_r  \rangle \, \bar{u}_r  \, \,dx \, dt | \ls
    \frac{1}{(N^2r)^{1-\sigma}} \frac{s^{-\frac{1}{4}}}{(Ns^{1/2})^{1-\sigma}} L^{\frac{7}{4}}
\end{equation}
for $ r \in (0,s) $ and the similar bound with permuted conjugation. By integration in $ r $, \eqref{eq:term1_bound} follows.
\begin{equation*}
    \text{LHS} \ \eqref{eq:term1_bound_reduced} \ls \|e^{(2s-r)\Delta} \Delta u \,  u_r\|_{L^2_{t,x}} \|\langle \nabla u_r \, , \nabla u_r  \rangle\|_{L^2_{t,x}}.
\end{equation*}
Then \eqref{eq:term1_bound_reduced} follows by combining  %\eqref{eq:term1_bd_reduced_1} and \eqref{eq:term1_bd_reduced_2} 
\begin{align} 
    \label{eq:term1_bd_reduced_1}
    \|\langle \nabla u_r, \nabla u_r  \rangle\|_{L^2_{t,x}} & \ls \vn{\nabla u_r}_{L^4_{t,x}}^2 \ls \frac{L}{(N^2r)^{1-\sigma}} ,\\
    \label{eq:term1_bd_reduced_2}
\|e^{(2s-r)\Delta} \Delta u \,  u_r\|_{L^2_{t,x}} & \ls 
 \frac{s^{-\frac{1}{4}}}{(Ns^{1/2})^{1-\sigma}} L^{\frac34}.
\end{align} 
The first bound \eqref{eq:term1_bd_reduced_1} follows from \eqref{eq:env_nabla_ur_L4}.
The bound \eqref{eq:term1_bd_reduced_2} is where the Bilinear Strichartz estimates from Corollary \ref{cor:bilinear_strichartz_variations} are used.
Based on $ \partial_s u_s = \Delta u_s $, split
\footnote{The significance of the threshold $ L^{-1}$ is that it corresponds to frequencies $ \ls 1 $ for the scaled solution (on $\mb{R}^2$).}
\[ u_r = u_{L^{-1}} - \int_r^{L^{-1}}  \Delta  u_l \, dl \]
The bound for $ \| \Delta u_{2s-r} \,  u_{L^{-1}}\|_{L^2_{t,x}} $ follows from 
\eqref{eq:inhom_bil}, \eqref{eq:envelope_bound_X1b}, \eqref{eq:envelope_bound_L4} and $ 2s-r \simeq s $. 

By splitting into cases $ 2s-r \leq l $ and $ 2s-r > l $, and applying \eqref{eq:inhom_bil} one has
\begin{equation}
    \vn{\Delta u_{2s-r} \, \Delta  u_l }_{L^2_{t,x}} \ls s^{-\frac{1}{4}} l^{-\frac{3}{4}} \vn{u_{s-r/2}}_{X^{1,b}} \vn{u_{l/2}}_{X^{1,b}}.
\end{equation} 
Using \eqref{eq:envelope_bound_X1b_all_s} for the last two factors, and integrating $ l $ over $ [r, L^{-1}] $ gives the second part of \eqref{eq:term1_bd_reduced_2}.

\subsection{Proof of the Term2 bound \eqref{eq:term2_bound}} \

Expressing $ e^{s\Delta}(|u|^2u) $ using \eqref{eq:commutator_def} and noting that the imaginary part of $ |\mathcal{D}_s(u)|^2 $ is $ 0$, one has
\begin{equation*} %\label{eq:term2_simplified}
\mathrm{Term}_{2,s} = - \Im \int_{t_0}^{t} \int_{\mathbb{H}^2} \overline{|u_s|^2 u_s}\, \mathcal{D}_s(u)\,dx\,dt,
\end{equation*}
which is bounded using H\"older's inequality by
\begin{equation} \label{eq:term2_holder}
|\mathrm{Term}_{2,s}| \ls \vn{u_s}_{L^\infty_{t,x}}^2 \, \vn{u_s}_{L^4_{t,x}} \, \vn{\mathcal{D}_s(u)}_{L^{4/3}_{t,x}} .
\end{equation}
To bound $ \mathcal{D}_s(u) $ use Lemma \ref{lem:commutator_structure}, \eqref{eq:envelope_bound_L4} and \eqref{eq:env_nabla_ur_L4}
\begin{equation*} %\label{eq:D_L43}
\vn{\mathcal{D}_s(u)}_{L^{4/3}_{t,x}} 
\ls  \int_0^s \vn{u_r}_{L^4_{t,x}} \vn{\nabla u_r}_{L^4_{t,x}}^2 \, dr 
\ls  L \int_0^s (N^2 r)^{-(1-\sigma)} \, dr 
\simeq L \, s \, (N^2 s)^{-(1-\sigma)} .
\end{equation*}
Combining this with \eqref{eq:envelope_bound_inf} and \eqref{eq:envelope_bound_L4} in \eqref{eq:term2_holder} concludes \eqref{eq:term2_bound}.

%% file: sections/proof_lwp.tex
\section{Proof of Lemmas and Propositions} \label{sec:proof_lemmas_prop}

\subsection{Proof of Proposition \ref{prop:modified_LWP} (Modified LWP)} \

Let $ u \in C([0,T_{\max}), H^{\sigma}) $ be the maximal $ H^{\sigma} $ solution given by the standard local theory. 
%Denote $ Q(t) := \vn{u}_{Y([0,t])} $ for $ t<T_{\max} $. 
By Lemma \ref{lem:Y_equiv}, it suffices to bound the equivalent norm $ \vn{u}_{\tilde{Y}([0,t])} $. 
The goal will be to prove that 
\begin{equation} \label{eq:bootstrap_LWP}
    \vn{u}_{\tilde{Y}([0,t])} \ls  L^{1/2} + \Gamma(t) \vn{u}_{\tilde{Y}([0,t])}^3, \qquad \Gamma(t) \leq t^{\alpha} L^{-\beta}, \qquad (1+\eta) \alpha+\beta \geq 1
%\Gamma(t) := t^{\sigma-3 \ep} L^{-\frac{2}{p}(1+\eta)+\frac{\eta}{4}}     
%Q(t) \ls  L^{1/2} + \Gamma(t) Q(t)^3
\end{equation}
for $ t \in J $.
Then, by a bootstrap argument for the continuous function $ q(t) := \vn{u}_{\tilde{Y}([0,t])}/L^{1/2} $ it follows that $ \vn{u}_{\tilde{Y}(J)} \ls L^{1/2} $ 
for $ J = [0, \delta L^{-1-\eta}] $ provided $ \delta $ is chosen small enough so that $ \Gamma \ll L^{-1} $.

\

It remains to prove \eqref{eq:bootstrap_LWP} for some appropriate exponents $ \alpha, \beta > 0 $. Choose $p,r$ by
\begin{equation} \label{eq:pr_numerology}
\frac{1}{p} = \frac{1-\sigma}{2} + \ep, \quad \frac{1}{r} = \frac{\sigma}{2} - \ep.
\end{equation}
Let $ u_s = e^{s\Delta} u $. Applying $ e^{s\Delta} $ to \eqref{eq:nls}, the equation for $ u_s $ can be written as
\begin{equation}
(i\partial_t + \Delta) u_s = |u_s|^2 u_s - \mathcal{D}_s(u), \qquad \mathcal{D}_s(u) := |u_s|^2 u_s-e^{s\Delta}(|u|^2 u).
\end{equation}

By \eqref{eq:X_embeddings} one has 
\begin{equation} \label{eq:F_emdedding}
    \vn{F}_{X^{0,b-1}} \ls \vn{F}_{L^{a} L^2}, \qquad  b=\frac{1}{2}+\ep, \quad a=(1-\ep)^{-1}=1+, \quad c=\big(\frac{1}{2}-\ep\big)^{-1}=2+.
\end{equation}
Therefore, by the inhomogeneous $X^{s,b}$ estimate \eqref{eq:linear_sol_estimate},
\begin{equation}
    \vn{u_s}_{X^{1,b}(I)} \ls \vn{u_s(0)}_{H^1} + \| |u_s|^2 u_s\|_{L^{a} H^1(I)} + \|\mathcal{D}_s(u)\|_{L^{a} H^1(I)}.
\end{equation}

Denote $ Q:= \vn{u}_{\tilde{Y}(I)} $, for a fixed interval $I$. The claim is 
\begin{align}
\| |u_s|^2 u_s\|_{L^{a} H^1(I)}  \ls \quad  \vn{u_s}_{X^{1,b}}  &  \vm{I}^{\alpha_1 }L^{-\beta_1} Q^2
\label{eq:nonlinearity_estimate_1} \\
\|\mathcal{D}_s(u)\|_{L^{a} H^1(I)}   \ls s^{\frac{1}{2}} (N^2 s)^{-(1-\sigma)} & \vm{I}^{\alpha_2} L^{-\beta_2+\frac{1}{2}} Q^3.
% s^{\frac{1}{2}} (N^2 s)^{-(1-\sigma)} \vm{I}^{\frac{1}{2}-} N^{0+} Q^3.
\label{eq:nonlinearity_estimate_2}
\end{align}
where
\[ \alpha_1 = \sigma-3 \ep, \quad \alpha_2 = \frac{\sigma}{2}- 2\ep, \qquad \beta_1 = \frac{2}{p}(1+\eta)-\frac{\eta}{4}, \quad \beta_2 = \frac{1}{2}+ \frac{1}{p}(1+\eta)-\frac{\eta}{8}. \
\]
Then \eqref{eq:bootstrap_LWP} follows by squaring these, multiplying with $ (N^2 s)^{1-\sigma} $ and integrating with $ ds/s$ and adding the corresponding bound for $ s_0 $. 
For the second bound, after integration, the $ s^{1/2} $ factor leads to a power $ N^{-1} $ and it suffices to use the crude bound $ N^{-1} \ls L^{-\frac{1}{2}}$. 

\begin{lemma} \label{lem:pr_embedding}
    Assume \eqref{eq:pr_numerology}, \eqref{eq:NLnumerology} hold. Then for every interval $I$ with $ \vm{I} \leq \delta L^{-1-\eta}$ one has 
\begin{equation} \label{eq:pinf_boundY}
\|u\|_{L_t^p L_x^\infty(I)} \lesssim L^{-\frac{1+\eta}{p}+\frac{\eta}{8}}  \vn{u}_{Y(I)}.
% L^{-\frac{\beta_1}{2}} \vn{u}_{Y(I)} \simeq L^{-\beta_2+\frac{1}{2}} \vn{u}_{Y(I)}
\end{equation}
% \[ \|u\|_{L_t^p L_x^\infty(I)} \lesssim \big( |I|^{1/p} N^\ep + N^{-(1-\sigma)} \big)  \vn{u}_{Y(I)}. \]
\end{lemma}
Note that the power of $ L $ equals $ -\frac{\beta_1}{2} = -\beta_2+\frac{1}{2} $.
By H\"older's inequality and Lemma \ref{lem:pr_embedding},
\begin{equation*}
\| |u_s|^2 u_s\|_{L^{a} H^1(I)}  \ls  \vm{I}^{\sigma-3 \ep} \vn{u_s}_{L^{\infty}H^1} \vn{u_s}_{L^p L^{\infty}}^2 
\ls \text{RHS} \ \eqref{eq:nonlinearity_estimate_1}
%\vm{I}^{\sigma-3 \ep} \vn{u_s}_{X^{1,b}} L^{-\frac{2}{p}(1+\eta)+\frac{\eta}{4}} Q^2. 
%\| |u_s|^2 u_s\|_{L^{a} H^1(I)}  \ls \vm{I}^{\sigma-3 \ep} \vn{u_s}_{L^{\infty}H^1} \vn{u_s}_{L^p L^{\infty}}^2 \ls \vm{I}^{\sigma-3 \ep} \big( |I|^{1/p} N^\ep + N^{-(1-\sigma)} \big)^2 Q^2. 
\end{equation*}

\

It remains to bound \eqref{eq:nonlinearity_estimate_2} for the commutator term $ \mathcal{D}_s(u) $, which is expressed by \eqref{eq:commutator_structure}. Fix $ r \in (0,s) $. 
Using \eqref{eq:strichartz_embedding} and the expansion \eqref{eq:heat_expansion}, similarly to \eqref{eq:envelope_bound_X1b}, one has
%the boundedness of the Riesz transform $ \nabla / (-\Delta)^{1/2} $ on $ L^q(\mathbb{H}^n) $ spaces, the Strichartz embedding $ X^{0,b} \subset L^4_{t,x} $ \eqref{} and then 
%splitting into spectral frequencies $ \leq N $ and $ \geq N $, it is seen that
\begin{equation} \label{eq:nabla_ur_L4}
\vn{\nabla u_r}_{L^4_{t,x}} \ls \vn{u_r}_{X^{1,b}} \ls (Nr^{\frac{1}{2}})^{-(1-\sigma)} Q.
\end{equation}

Applying the following bound, based on \eqref{eq:F_emdedding},
\begin{equation} \label{eq:F_emdedding_heat}
\vn{e^{(s-r)\Delta}F}_{X^{1,b-1}} \ls \vn{e^{(s-r)\Delta}F}_{L^{a} H^1} \ls (s-r)^{-\frac{1}{2}} \vn{F}_{L^{a} L^2}.
\end{equation}
with $ F= \langle \nabla u_r, \nabla \bar{u}_r \rangle u_r $ and $ \langle \nabla u_r, \nabla  u_r  \rangle \bar{u}_r $ implies,
using \eqref{eq:nabla_ur_L4}, H\"older and Lemma \ref{lem:pr_embedding}
\begin{equation*}
\eqref{eq:F_emdedding_heat} \ls (s-r)^{-\frac{1}{2}} \vn{\nabla u_r}_{L^4_{t,x}}^2 \vn{u_r}_{L^{c} L^\infty} 
\ls (s-r)^{-\frac{1}{2}} (N^2r)^{-(1-\sigma)} \vm{I}^{\alpha_2} L^{-\beta_2+\frac{1}{2}} Q^3
\end{equation*}
Then \eqref{eq:nonlinearity_estimate_2} follows from the basic estimate
\[ \int_0^s  (s-r)^{-\frac{1}{2}} (N^2r)^{-(1-\sigma)} \, dr \, \simeq \, s^{\frac{1}{2}} (N^2 s)^{-(1-\sigma)}.
\]
To conclude the first part, since $ Y([0,t]) $ controls $ \vn{u}_{L^{\infty} H^{\sigma}([0,t])} $, the bound on $ [0,t]$ and the standard continuation criterion imply $ T_{\max} > \delta  L^{-1-\eta}$.

\

It remains to establish the bound $ \|u\|_{X^{0,b}(J)}\lesssim_{M} 1 $ for $ M=\|u_0\|_{L^2} $. 
By \eqref{eq:linear_sol_estimate} and \eqref{eq:F_emdedding}
\[ \|u\|_{X^{0,b}(J)} \lesssim \|u_0\|_{L^2} + \||u|^2u\|_{L^aL^2(J)} 
\ls \|u_0\|_{L^2} + \|u\|_{L_t^\infty L_x^2(J)}\|u\|_{L_t^{2a}L_x^\infty(J)}^2
\]
Hence, by mass conservation, 
$ \|u\|_{X^{0,b}(J)} \lesssim M+M\|u\|_{L_t^{2a}L_x^\infty(J)}^2 $. 
By Hölder in time, Lemma \ref{lem:pr_embedding}, and $  \vn{u}_{Y(I)} \ls L^{\frac{1}{2}} $ the numerology shows that 
$ \|u\|_{L_t^{2a}L_x^\infty(J)} $ is bounded by a negative power of $ L $, 
concluding $  \|u\|_{X^{0,b}(J)}\lesssim_M 1 $.

%%%%%%%%%%%%%%%%%%%%%%%%%%%%%%%%%%%%%%%%%%%%%%%%%%%%%%%
%%%%%%%%%%%%%%%%%%%%%%%%%%%%%%%%%%%%%%%%%%%%%%%%%%%%%%%

\

\subsection{Proof of Lemma \ref{lem:pr_embedding}} \ 

Split $ u = u_{\mathrm{lo}} + u_{\mathrm{hi}} $ where $ u_{\mathrm{lo}}, u_{\mathrm{hi}} $ contain spectral frequencies $ \ls N $ and $ \gg N $ respectively. 
For the low piece, use $H^{1+\ep} \subset L^\infty$, H\"older in $ t $ and the embedding $ X^{1,\frac12+} \subset L^\infty H^1 $ 
\begin{equation*}
\|u_{\mathrm{lo}}\|_{L_t^p L_x^\infty(I)} \le |I|^{1/p} \|u_{\mathrm{lo}}\|_{L_t^\infty H^{1+\ep}(I)} \lesssim |I|^{1/p} N^\ep \|u_{\mathrm{lo}}\|_{X^{1,b}(I)} \ls |I|^{1/p} N^\ep \vn{e^{N^{-2} \Delta} u}_{X^{1,b}(I)}.
\end{equation*}

For the high piece, use $W^{\sigma,r} \subset L^\infty$ and the Strichartz embedding $ X^{\sigma,b} \subset L^p W^{\sigma,r} $
\begin{equation*}
\|u_{\mathrm{hi}}\|_{L_t^p L_x^\infty(I)} \lesssim \|u_{\mathrm{hi}}\|_{L_t^p W_x^{\sigma,r}(I)} \lesssim \|u_{\mathrm{hi}}\|_{X^{\sigma,b}(I)} \lesssim N^{-(1-\sigma)} \vn{u}_{Y(I)}.
\end{equation*}
Due to \eqref{eq:NLnumerology} both terms are dominated by the RHS of \eqref{eq:pinf_boundY}.
$\square$

%%%%%%%%%%%%%%%%%%%%%%%%%%%%%%%%%%%%%%%%%%%%%%%%%%%%%%%%
%%%%%%%%%%%%%%%%%%%%%%%%%%%%%%%%%%%%%%%%%%%%%%%%%%%%%%%%

\subsection{Proof of Lemma \ref{lem:energy_comparison}}

The first part of \eqref{eq:energy_comparison_1} follows from the triangle inequality and $ H^1\subset H^\sigma$. 
For the second part, first note that $ P u(t)= u_{s_0}(t)$ and its $ H^1 $ norm is controlled by $E(t,s_0)$, which is part of $ \mathcal{IE}_N(t)$.
To control $ (1-P) u(t) $ by $ \mathcal{IE}_N(t) $, use the bound 
\[ %\label{eq:high_symbol_bound}
N^{2(1-\sigma)}x^\sigma(1-e^{-x/N^2})^2
\lesssim
\int_0^{N^{-2}} (N^2s)^{1-\sigma} x e^{-2sx}\,\frac{ds}{s},
\]
where $ x $ denotes the spectral variable of $ -\Delta_{\mathbb{H}^2} $.

To prove \eqref{eq:energy_comparison_2}, first note that the $ \vm{u_s}^4 $ part of $ \mathcal{IE}_N(t) $ follows from the boundedness of $ e^{s\Delta} $ on $ L^4 $.
The gradient part of \eqref{eq:energy_comparison_2} follows from the bound
\[
\sup_{0<s\le s_0}(N^2s)^{1-\sigma}x e^{-2sx} +
\int_0^{s_0}(N^2s)^{1-\sigma}x e^{-2sx}\,\frac{ds}{s}
\lesssim N^{2(1-\sigma)}x^\sigma. \qquad \square
\]

%% file: sections/appendix.tex
\section{Appendix} 
%\section{Harmonic analysis and Function spaces} \

\subsection{The Helgason-Fourier transform}

Let $x \in \mathbb{H}^n$, $\omega \in S^{n-1}$ (the boundary), and $\lambda \in \mathbb{R}$ (the spectral parameter). 
Let $\langle x, \omega \rangle$ denote the horocycle bracket, i.e. the signed distance from the origin to the horosphere passing through $x$ determined by $\omega$. 
Denote $\rho = \frac{n-1}{2}$.

The Helgason-Fourier transform \cite{Helgason1, helgason2001differential, helgason2024geometric} $ \mathcal{F}_H f $ of a function $f \in C_c^\infty(\mathbb{H}^n)$ is defined as:
\begin{equation} \label{eq:helgason_fourier_def}
    \widetilde{f}(\lambda, \omega) = \int_{\mathbb{H}^n} f(x) \, e^{(-i\lambda + \rho)\langle x, \omega \rangle} \, dx
\end{equation}
The inversion formula expresses $f$ as a superposition of eigenfunctions:
\begin{equation} \label{eq:inversion}
    f(x) = \int_{\mathbb{R_{+}}} \int_{\partial \mathbb{H}^n} \widetilde{f}(\lambda, \omega) e^{(i\lambda + \rho)\langle x, \omega \rangle} |c(\lambda)|^{-2} d\omega d\lambda,
\end{equation}
where $ d \omega $ is the $ K=SO(n)$ -invariant probability measure on $ \partial \mathbb{H}^n $ and $ c(\lambda) $ is the Harish-Chandra $ c $-function 
\begin{equation}\label{eq:c:function}
c(\lambda) = C \frac{\Gamma(i\lambda)}{\Gamma(i\lambda+\rho)},\qquad
|c(\lambda)|^{-2} = C \frac{|\Gamma(i\lambda+\rho)|^{2}}{|\Gamma(i\lambda)|^{2}} \asymp \lmd^2 (1+|\lmd|)^{n-3} .
\end{equation}
The Plancherel formula extends $ \mathcal{F}_H f $ to a unitary isomorphism between Hilbert spaces:
\begin{equation}\label{eq:plancherel}
    \int_{\mathbb{H}^n} |f(x)|^2 \, dx =  \int_{\mathbb{R}_+} \int_{S^{n-1}} |\widetilde{f}(\lambda, \omega)|^2 \, \frac{d\omega \, d\lambda}{|c(\lambda)|^2}
\end{equation}
%Alternatively, one can restrict to the positive spectrum $\mathbb{R}_+$ and remove the $ 1/2$ in \eqref{eq:plancherel}. 

Now consider functions in $L^2(\mathbb{R} \times \mathbb{H}^n)$. One defines the combined transform $\mathcal{F} = \mathcal{F}_t \otimes \mathcal{F}_H$ 
\begin{equation} \label{eq:combined_transform}
    \mathcal{F} f(\tau, \lambda, \omega) = \int_{-\infty}^{\infty}   e^{-i\tau t} \, \int_{\mathbb{H}^n} f(t, x) \, e^{(-i\lambda + \rho)\langle x, \omega \rangle} \, dx dt
\end{equation}
One may take $ \mathcal{F}_t $ as the $L^2(\mathbb{H}^n)$-valued Fourier transform (Bochner integral) of the function $ f(t) $. 
Then by Parseval's identity for vector-valued functions, the map $\mathcal{F}$ becomes a unitary operator
\begin{equation}\label{eq:FT}
    \mathcal{F} : L^2(\mathbb{R}_t \times \mathbb{H}^n_x) \longrightarrow L^2\left(\mathbb{R}_\tau \times \mathbb{R}^+_\lambda \times S^{n-1}_\omega, \, |c(\lambda)|^{-2} \, d\tau \, d\omega \, d\lambda\right)
\end{equation}

\subsection{Spherical analysis} \

A function $ f $ on $ \mathbb{H}^n $ is called radial if it depends only on the geodesic distance from the origin, equivalently if 
$ f(k \cdot x) = f(x) $ for all $ k \in K= SO(n) $. Identifying $ \mathbb{H}^n = G/K $ as a homogeneous space 
where $ G=SO_0(n,1) $, radial functions are lifted to $ K $-bi-invariant functions on $ G $.

Based on the group convolution on $ G $, it becomes meaningful to define the convolution of functions $ f $ on $ \mathbb{H}^n $ 
only with a radial kernel $ k(x)=\kappa(d(o,x)) $ where $\kappa: [0,\infty) \to \mathbb{C} $, as:
\begin{equation} \label{eq:convolution_def}
    T_kf(x) = (f*k)(x) := \int_{G} f(g \cdot o) k(g^{-1} \cdot x)  \, dg = \int_{\mathbb{H}^n} f(y) \kappa(d(x,y)) dy.
\end{equation}
The spherical functions $ \varphi_\lambda $ on $ \mathbb{H}^n $ are defined as normalized smooth radial eigenfunctions of $ \Delta_{\mathbb{H}^n} $ 
\begin{equation} \label{eq:spherical_functions_1}
    \Delta  \varphi_\lambda = - (\lambda^2 + \rho^2) \varphi_\lambda, \qquad \varphi_\lambda(0) = 1, \qquad \quad \lambda \in \mathbb{C}, \quad \rho = \frac{n-1}{2}. 
\end{equation}
These are obtained by averaging the Helgason waves over the boundary
\begin{equation} \label{eq:spherical_functions_2}
    \varphi_\lambda(x) = \int_{\partial \mathbb{H}^n} e^{(i\lambda + \rho)\langle x, \omega \rangle} d\omega, \qquad \varphi_\lambda = \varphi_{-\lambda}.
\end{equation}
For $ n=3 $ one has 
\begin{equation}\label{eq:3d:spherical}
\varphi_\lambda(r)=\frac{\sin(\lambda r)}{\lambda\,\sinh r}, \qquad \text{on} \quad \mathbb{H}^3.
\end{equation}
More generally, the spherical functions can be expressed in terms of the Jacobi functions $ \phi_{\lambda}^{\alpha,\beta} $  and  
\cite[{(15.9.11)}]{NIST:DLMF} the Gauss hypergeometric function ${}_2F_1$ as     
\begin{equation}\label{eq:Jacobi:Gauss}
\varphi_\lambda(r) = \phi_{\lambda}^{\rho-\frac{1}{2},-\frac{1}{2}} (r) = 
{}_2F_1\!\left(\frac{\rho + i \lambda}{2},\,
\frac{\rho - i\lambda}{2};\,
\frac{n}{2};\,
-\sinh^2 r\right),
\end{equation}
If $ k $ is a radial function, the Helgason-Fourier transform \eqref{eq:helgason_fourier_def} and the inversion \eqref{eq:inversion} simplify as
\begin{equation} \label{eq:spherical_transform}
    \widetilde{k}(\lambda, \omega) = \widetilde{k}(\lambda) =  \int_{\mathbb{H}^n} k(x) \, \varphi_{-\lambda}(x) \, dx, \qquad     k(x) = \int_{\mathbb{R}_+} \widetilde{k}(\lambda) \, \varphi_{\lambda}(x) \, |c(\lambda)|^{-2} \, d\lambda.
\end{equation}
Additionally, the convolution \eqref{eq:convolution_def} transforms as  
\begin{equation} \label{eq:transformed_convolution}
  \widetilde{T_kf}(\lambda,\omega) = \widetilde{f*k}(\lambda,\omega) = \widetilde{f}(\lambda, \omega) \, \widetilde{k}(\lambda), 
  %\quad \text{and} \quad     \|T_k\|_{L^2(\mathbb H^n) \to L^2(\mathbb H^n)} = \|\widetilde{k}\|_{L^\infty(\mathbb R_+)} .
\end{equation}
In particular, spectral multipliers such as $ e^{-it\Delta}, e^{s\Delta}, P_N=\chi(-\Delta/N^2) $ are defined in this way. 
Similarly, a space-time integral operator $ T_k $ with kernel depending only on $ t-s $ and $ d(x,y) $ is unitarily equivalent under \eqref{eq:FT} to multiplication 
%by a function $ m(\tau, \lambda) $
\begin{equation} \label{eq:space_time_multiplier}
    T_k =  \mathcal{F}^{-1} m(\tau, \lambda) \mathcal{F}, \qquad \text{and} \qquad \|T_k\|_{L^2(\mb{R} \times \mathbb H^n) \to L^2(\mb{R} \times \mathbb H^n)} = \|m\|_{L^\infty(\mb{R} \times \mathbb R_+)}.
\end{equation}

\subsection{The $X^{s,b}$ spaces} \label{sec:xsb_spaces}
These hyperbolic space analogues of the Bourgain spaces \cite{bourgain1993fourier} will be defined by simply replacing the euclidean Fourier transform with the Helgason-Fourier transform \eqref{eq:combined_transform}. \\
For a general approach to $X^{s,b}$ on curved backgrounds, see \cite{gavrus2021wave}, \cite{tataru1996xs}.
% A more general treatment of $ X^{s,b} $ spaces on curved backgrounds, not based on symmetries but rather on interpolation, appears in \cite{}, \cite{}. 
On compact manifolds see \cite{burq2005bilinear}.

Let $ s, b \in \mathbb{R} $. The $X^{s,b}$ spaces on $\mathbb{H}^n$ are defined by the norm 
$$
\|u\|_{X^{s,b}(\mathbb{H}^n)} := \left\| \jb{\lmd}_{\rho}^s \langle  \tau + \jb{\lmd}^2_{\rho} \rangle^b \widetilde{u}(\tau, \lambda, \omega) \right\|_{L^2_{\tau,\lambda,\omega}(d\tau d\mu)}.
$$
where $ d\tau d\mu $ is the measure from \eqref{eq:FT} and $ \jb{\lmd}_{\rho} = \sqrt{\lambda^2 + \rho^2}$ with $ \jb{\lmd}_{\rho} \simeq \langle \lambda \rangle$. 
This means that $ \jb{\lmd}_{\rho} $ and $ \jb{\lmd} $ are interchangeable. Given an interval $ I $ one identifies two functions $ u \sim v $ if they agree on $ I $. 
One defines the restriction space $ X_I^{s,b}(\mathbb{H}^n) $ as $ X^{s,b}(\mathbb{H}^n) / \sim$, endowed with the norm
$$ \vn{u}_{X_I^{s,b}(\mathbb{H}^n)} := \inf \{ \vn{\tilde{u}}_{X^{s,b}(\mathbb{H}^n)} : \tilde{u} \in X^{s,b}(\mathbb{H}^n), \, \tilde{u} \sim u \}.$$
The space $X_I^{s,b}(\mathbb{H}^n)$ is complete and embeds into $C(I, H^s(\mathbb{H}^n)) $ when $ b > \frac12 $, as follows. 

\begin{proposition}\label{prop:X_properties} Let $ s, b \in \mathbb{R} $. Then the following properties hold:
\begin{enumerate}
    \item (Isometry) Defining $  v(t) := e^{- it\Delta} u(t) $ one has the isometry with product Sobolev spaces
    \begin{equation}\label{eq:isometry_product}
        \|u\|_{X^{s,b}(\mathbb{H}^n)} = \|v\|_{H^b_t H^s_x(\mathbb{R} \times \mathbb{H}^n)}
    \end{equation}
    \item (Duality) With the $L^2(\mathbb{R} \times \mathbb{H}^n)$ inner product pairing one has  
    \begin{equation} \label{eq:duality}
    |\langle u, v \rangle| \le \|u\|_{X^{s,b}(\mathbb{H}^n)} \|v\|_{X^{-s,-b}(\mathbb{H}^n)}, \qquad  X^{s,b}(\mathbb{H}^n) ^* \simeq X^{-s,-b}(\mathbb{H}^n).
    \end{equation}
    \item (Interpolation) The intermediate spaces are recovered by complex interpolation 
    \begin{equation}\label{eq:complex_interp}
        X^{s,\theta}(\mathbb{H}^n) = [X^{s,0}(\mathbb{H}^n), X^{s,1}(\mathbb{H}^n)]_\theta, \qquad \theta \in (0, 1).
    \end{equation}
    \item (Time localization) If $\eta \in C_c(\mathbb{R})$ then $ \|\eta(t)u\|_{X^{s,b}(\mathbb{H}^n)} \lesssim \|u\|_{X^{s,b}(\mathbb{H}^n)} $. Moreover, for $ b>1/2$ 
    \begin{equation} \label{eq:linear_sol_estimate}
        (Linear \ estimate) \qquad \quad \ \,
    \|\eta(t)u\|_{X^{s,b}(\mathbb{H}^n)} \lesssim \vn{u(0)}_{H^s} + \|(i\pt_t+\Delta)u\|_{X^{s,b-1}(\mathbb{H}^n)}.
    \end{equation}
    \item (Decomposition into free solutions) Let $u \in X^{s,b}(\mathbb{H}^n)$ for $b > 1/2$. 
    Then there exist $f_{h} \in H^s(\mathbb{H}^n)$, $ h \in \mathbb{R}$ such that 
    \begin{equation}\label{eq:integral_rep}
    u(t) = \int_{\mathbb{R}} e^{ith} e^{it\Delta} f_{h} \dd h, \quad  \text{with}  \quad
    \int_{\mathbb{R}} \|f_{h}\|_{H^s(\mathbb{H}^n)} \dd h %\lesssim  \int_{\mathbb{R}} \langle h \rangle^{2b} \|f_{h}\|_{H^s(\mathbb{H}^n)}^2 \dd h ^{1/2} 
    \ls \|u\|_{X^{s,b}(\mathbb{H}^n)},
    \end{equation}
     \item (Embeddings) Consequently, $ X^{s,b}(\mathbb{H}^n) \subset C(\mathbb{R}, H^s(\mathbb{H}^n)) $ and for $ b = \frac{1}{2}+\ep $, $ a=(1-\ep)^{-1} $ 
    % $ \|u\|_{L^\infty_t H^s} \lesssim \|u\|_{X^{s,b}}. $
    \begin{equation} \label{eq:X_embeddings}
     \|u\|_{L^\infty_t H^s(\mathbb{H}^n)} \lesssim \|u\|_{X^{s,b}(\mathbb{H}^n)}, \qquad
     \vn{F}_{X^{0,b-1}} \ls \vn{F}_{L^{a} L^2}
    \end{equation}
    \item (Strichartz embedding) For $ b > 1/2 $, $ p \geq 2, \ q>2 $, one has
    \begin{equation} \label{eq:strichartz_embedding}
    \vn{u}_{L^p_t L^q_x(\mathbb{R} \times \mathbb{H}^n)} \lesssim \|u\|_{X^{0,b}(\mathbb{H}^n)}, \quad \vn{\nabla u}_{L^p_t L^q_x(\mathbb{R} \times \mathbb{H}^n)} \lesssim \|u\|_{X^{1,b}(\mathbb{H}^n)},   \quad \frac{2}{p} \geq \frac{n}{2}-\frac{n}{q}.
    \end{equation}
\end{enumerate}
\end{proposition}

\begin{proof}
These are proved by following the corresponding euclidean versions and replacing $ |\xi| $ by $ \jb{\lmd}_{\rho} $. 
\begin{enumerate}
    \item Note that $  \tilde{v}(\tau, \lambda, \omega) = \tilde{u}(\tau - \jb{\lmd}^2_\rho, \lambda, \omega) $. Then \eqref{eq:isometry_product} follows by a change of variables. 
    \item Starting with $ |\langle u, v \rangle|$, use Plancherel's theorem \eqref{eq:FT}, 
distribute the weights $\jb{\lmd}_\rho^s \langle \tau + \jb{\lmd}_\rho^2 \rangle^b$ and their inverses
and apply the Cauchy-Schwarz inequality.
    \item They are weighted $L^2$ spaces with $w(s,\theta) = \jb{\lmd}_\rho^s \langle \tau + \jb{\lmd}_\rho^2 \rangle^\theta$ and 
$w(s,\theta) = w(s,0)^{1-\theta} w(s,1)^\theta$. 
    \item See Lemma 2.11 and Prop. 2.12 in \cite{tao2006nonlinear}.
    \item Express $\tilde{u}(\tau, \lambda, \omega) = \int \tilde{u}(-\jb{\lmd}_\rho^2 + h, \lambda, \omega) \delta(\tau + \jb{\lmd}_\rho^2 - h ) \dd h. $
For each $h \in \mathbb{R}$, one defines the spatial function $f_{h} \in H^s(\mathbb{H}^n)$ by
$ \tilde{f}_{h}(\lambda, \omega) := \tilde{u}(-\jb{\lmd}_\rho^2 + h, \lambda, \omega).$
The space-time transform of $ e^{ith} e^{it\Delta} f_{h}$ is exactly $\tilde{f}_{h}(\lambda, \omega) \delta(\tau+ \jb{\lmd}_\rho^2- h )$. 

To establish the bounds, change variables to $\tau = - \jb{\lmd}_\rho^2 + h$:
\begin{equation*}
    \|u\|_{X^{s,b}(\mathbb{H}^n)}^2 
    = \int_{\mathbb{R}} \langle h \rangle^{2b} \Big( \int_{\mathbb{R}^+ \times \mathbb{S}^{n-1}} \jb{\lmd}_\rho^{2s} |\tilde{f}_{h}(\lambda, \omega)|^2 \dd\mu(\lambda, \omega) \Big) \dd h
    = \int_{\mathbb{R}} \langle h \rangle^{2b} \|f_{h}\|_{H^s(\mathbb{H}^n)}^2 \dd h.
\end{equation*}
Since $b > 1/2$, an application of the Cauchy-Schwarz inequality in $h$ yields \eqref{eq:integral_rep}.
% \begin{equation*}
%     \int_{\mathbb{R}} \|f_{h}\|_{H^s(\mathbb{H}^n)} \dd h \le \left( \int_{\mathbb{R}} \langle h \rangle^{-2b} \dd h \right)^{1/2} \left( \int_{\mathbb{R}} \langle h \rangle^{2b} \|f_{h}\|_{H^s(\mathbb{H}^n)}^2 \dd h \right)^{1/2} \lesssim \|u\|_{X^{s,b}(\mathbb{H}^n)}.
% \end{equation*}
    \item The first bound follows from \eqref{eq:integral_rep}, the unitarity of $e^{it(h+\Delta)}$ on $H^s(\mathbb{H}^n)$ and from Minkowski's inequality in $h$.
For the continuity $ C(\mathbb{R}, H^s(\mathbb{H}^n))$ use the dominated convergence theorem. \\
To prove the second bound, use the isometry from (1) together with $ L^{a}_t \subset H^{-1/2+\varepsilon}_t $, 
which is the dual of the Sobolev embedding $ H^{\frac12-\varepsilon}(\mb{R}) \subset L^{1/\ep}(\mb{R}) $,
and use the unitarity of the flow.

\item In \cite{anker2009nonlinear,banica2007nonlinear,ionescu2009semilinear} Strichartz estimates % \cite{anker2009nonlinear}, \cite{banica2007nonlinear}, \cite{ionescu2009semilinear}
$ e^{it\Delta}: L^2_x(\mathbb{H}^n) \to L^p_t L^q_x(\mathbb{R} \times \mathbb{H}^n) $ are proved in the triangle range $ \frac{2}{p} \geq \frac{n}{2}-\frac{n}{q}, \ p \geq 2, \ q>2 $. 
The first bound in \eqref{eq:strichartz_embedding} follows from applying these for each $ f_h $ in the representation \eqref{eq:integral_rep} and by Minkowski's inequality in $ h $. 
For the second bound use the boundedness of the Riesz transform $ \nabla / (-\Delta)^{1/2} $ on $ L^q(\mathbb{H}^n) $ spaces \cite[Theorem 6.1]{strichartz1983analysis}. 

% \begin{equation}
% \vn{e^{it\Delta} f}_{L^p_t L^q_x(\mathbb{R} \times \mathbb{H}^n)} \lesssim \|f\|_{L^2_x(\mathbb{H}^n)}, \qquad \qquad \frac{2}{p} \geq \frac{n}{2}-\frac{n}{q}.
% \end{equation}
\end{enumerate}
\end{proof}